\documentclass[12pt]{amsart}      %{article} was 12pt latex e
\usepackage{txfonts}      %{article} was 12pt latex e

\usepackage{amssymb}
\usepackage{eucal}
\usepackage{amsmath}
\usepackage{amscd}
\usepackage{multicol}
\usepackage[all]{xy}           %xypic macro for latex2.09
\usepackage{graphicx}
\usepackage{color}
\usepackage{colordvi}
\usepackage{xspace}
\usepackage{epsfig}
\usepackage{float}
\usepackage{hyperref}

\newcommand{\nc}{\newcommand}
\newcommand{\delete}[1]{}

\nc{\dfootnote}[1]{{}}          %{{}}
\nc{\ffootnote}[1]{\dfootnote{#1}}

\delete{
\nc{\mfootnote}[1]{{}}        % Use this to suppress footnotes
\nc{\ofootnote}[1]{{}}        % Use this to suppress footnotes
}

\nc{\mfootnote}[1]{\footnote{#1}} % Use this to show footnotes
\nc{\ofootnote}[1]{\footnote{\tiny Older version: #1}} % Use this to show footnotes

\nc{\mlabel}[1]{\label{#1}}  % Use this to suppress names
\nc{\mcite}[1]{\cite{#1}}  % Use this to suppress names
\nc{\mref}[1]{\ref{#1}}  % Use this to suppress names
\nc{\mkeep}[1]{{}}      % Use this to suppress marginpar
\nc{\mbibitem}[1]{\bibitem{#1}} % Use this to show number name
\delete{
\nc{\mcite}[1]{\cite{#1}{{\bf{{\ }(#1)}}}}  % Use this lines to show names
\nc{\mlabel}[1]{\label{#1}  % Use the next two lines to show names
{\hfill \hspace{1cm}{\bf{{\ }\hfill(#1)}}}}
\nc{\mref}[1]{\ref{#1}{{\bf{{\ }(#1)}}}}  % Use this lines to show names
\nc{\mbibitem}[1]{\bibitem[\bf #1]{#1}} % Use this to show name
\nc{\mkeep}[1]{\marginpar{{\bf #1}}} % Use this to show marginpar
}

\nc{\qssha}{{{\ssha\hspace{-2pt}_\ast}}\,}
\newtheorem{theorem}{Theorem}[section]
\newtheorem{prop}[theorem]{Proposition}
\newtheorem{defn}[theorem]{Definition}

\newtheorem{coro}[theorem]{Corollary}
\newtheorem{prop-def}{Proposition-Definition}[section]

\newtheorem{exam}[theorem]{Example}

\newtheorem{notation}[theorem]{Notation}

\nc{\bond}{\vdash}

\nc{\un}{u}                 %unit map in bialgebra
\nc{\mult}{m}       %multiplication in bialgebra
\nc{\cprod}{\star}

\nc{\comp}[1]{\langle #1\rangle}
\nc{\spr}{\cdot}
\nc{\disp}[1]{\displaystyle{#1}}
\nc{\bin}[2]{ (_{\stackrel{\scs{#1}}{\scs{#2}}})}  %binomial coeff
\nc{\binc}[2]{ \left (\!\! \begin{array}{c} \scs{#1}\\
    \scs{#2} \end{array}\!\! \right )}  %binomial coeff
\nc{\bbinc}[2]{ \left (\!\! \begin{array}{c} {#1}\\
    {#2} \end{array}\!\! \right )}  %binomial coeff
\nc{\bincc}[2]{  \left ( {\scs{#1} \atop
    \vspace{-.5cm}\scs{#2}} \right )}  %binomial coeff
\nc{\sarray}[2]{\begin{array}{c}#1 \vspace{.1cm}\\ \hline
    \vspace{-.35cm} \\ #2 \end{array}}
\nc{\bs}{\bar{S}} \nc{\dcup}{\stackrel{\bullet}{\cup}}
\nc{\dbigcup}{\stackrel{\bullet}{\bigcup}} \nc{\grd}[1]{^{(#1)}}
\nc{\mzeta}{\zeta} \nc{\PF}{\mathbf{PF}} \nc{\la}{\longrightarrow}
\nc{\fe}{\'{e}} \nc{\rar}{\rightarrow} \nc{\dar}{\downarrow}
\nc{\dap}[1]{\downarrow \rlap{$\scriptstyle{#1}$}}
\nc{\uap}[1]{\uparrow \rlap{$\scriptstyle{#1}$}}
\nc{\defeq}{\stackrel{\rm def}{=}} \nc{\dis}[1]{\displaystyle{#1}}
\nc{\dll}{_{\mrm{dll}}} \nc{\dotcup}{\
\displaystyle{\bigcup^\bullet}\ } \nc{\hcm}{\ \hat{,}\ }
\nc{\hcirc}{\hat{\circ}} \nc{\hts}{\hat{\shpr}}
\nc{\lts}{\stackrel{\leftarrow}{\shpr}}
\nc{\rts}{\stackrel{\rightarrow}{\shpr}} \nc{\lleft}{[}
\nc{\lright}{]} \nc{\uni}[1]{\tilde{#1}} \nc{\free}[1]{\bar{#1}}
\nc{\gzeta}{\bar{\zeta}}
\nc{\den}[1]{\check{#1}} \nc{\lrpa}{\wr} \nc{\curlyl}{\left \{
\begin{array}{c} {} \\ {} \end{array}
    \right . \!\!\!\!\!\!\!}
\nc{\curlyr}{ \!\!\!\!\!\!\!
    \left . \begin{array}{c} {} \\ {} \end{array}
    \right \} }
\nc{\longmid}{\left | \begin{array}{c} {} \\ {} \end{array}
    \right . \!\!\!\!\!\!\!}
\nc{\ot}{\otimes} \nc{\bigot}{\bigotimes} \nc{\mdiv}{\mrm{div}}
\nc{\qspr}{\ast}
\nc{\sg}{S} \nc{\ig}{I} \nc{\pg}{P} \nc{\jg}{J} \nc{\eg}{E}
\nc{\fg}{F} \nc{\cg}{C} \nc{\mg}{M} \nc{\abg}{C} \nc{\bas}{B}
\nc{\Lyn}{\mrm{Lyn}} \nc{\Li}{\mrm{Li}}
\nc{\lyn}{\Lyn}
\nc{\sG}{{\cals}} \nc{\iG}{{\cali}} \nc{\jG}{{\calj}}
\nc{\eG}{{\cale}} \nc{\pG}{{\calp}} \nc{\fG}{{\calf}}
\nc{\cG}{{\calc}} \nc{\mG}{{\calm}}
\nc{\ora}[1]{\stackrel{#1}{\rar}}
\nc{\ola}[1]{\stackrel{#1}{\la}}%${\Bbb Z}$
\nc{\pex}[1]{\{#1\}}
\nc{\scs}[1]{\scriptstyle{#1}}
\nc{\mrm}[1]{{\rm #1}}
\nc{\sym}[1]{{\widehat{#1}}}
\nc{\margin}[1]{\marginpar{\rm #1}}   %{\rm #1}}
\nc{\dirlim}{\displaystyle{\lim_{\longrightarrow}}\,}
\nc{\invlim}{\displaystyle{\lim_{\longleftarrow}}\,}
\nc{\mvp}{\vspace{0.5cm}}
\nc{\svp}{\vspace{2cm}}
\nc{\vp}{\vspace{8cm}}
\nc{\proofbegin}{\noindent{\bf Proof: }}
\nc{\proofend}{$\blacksquare$ \vspace{0.5cm}}
\nc{\shqs}{\eta}
\font\cyr=wncyr10
\nc{\sha}{{\mbox{\cyr X}}}  %used to be \cyr
\newfont{\scyr}{wncyr10 scaled 550}
\nc{\ssha}{\mbox{\bf \scyr X}}
\nc{\lzero}{_{\hskip -5pt
0}} \nc{\shzero}{_{\hskip -7.5pt 0}}
\newfont{\bcyr}{wncyr10 scaled 1000}
\nc{\ncsha}{{\mbox{\cyr X}^{\mathrm NC}}}
\nc{\ncshao}{{\mbox{\cyr X}^{\mathrm NC,\,0}}}
\nc{\shpr}{\diamond}    %Shuffle product

\nc{\shf}{{^{\ssha}}}
\nc{\qsh}{{^{\ast}}}
\nc{\esh}{{^{ext}}}
\nc{\lshf}{_{\ssha}}
\nc{\lqsh}{_{\ast}}
\nc{\shprl}{{{\shpr}_\lambda}}
\nc{\shpro}{\diamond^0}    %Shuffle product
\nc{\shpru}{\check{\diamond}}
\nc{\catpr}{\diamond_l}
\nc{\rcatpr}{\diamond_r}
\nc{\lapr}{\diamond_a}
\nc{\lepr}{\diamond_e}
\nc{\tcon}{^{\ot}}
\nc{\conv}{_c}
\nc{\vep}{\varepsilon}
\nc{\labs}{\mid\!}
\nc{\rabs}{\!\mid}
\nc{\hsha}{\widehat{\sha}}
\nc{\lsha}{\stackrel{\leftarrow}{\sha}}
\nc{\rsha}{\stackrel{\rightarrow}{\sha}}
\nc{\EDS}{{\mathrm{EDS}}}
\nc{\DS}{{\mathbf{DS}}}
\nc{\lc}{[}
\nc{\rc}{]}
\nc{\rbset}{R}
\nc{\rbnum}{r}
\nc{\rbfun}{\mathbf{R}}
\nc{\pset}{P}
\nc{\pnum}{p}
\nc{\pfun}{\mathbf{P}}
\nc{\spset}{SP}
\nc{\spnum}{sp}
\nc{\spgen}{\mathbf{SP}}
\nc{\srbi}[1]{\{#1\}}

\nc{\ann}{\mrm{ann}}
\nc{\Aut}{\mrm{Aut}}
\nc{\can}{\mrm{can}}
\nc{\colim}{\mrm{colim}}
\nc{\Cont}{\mrm{Cont}}
\nc{\rchar}{\mrm{char}}
\nc{\cok}{\mrm{coker}}
\nc{\dtf}{{R-{\rm tf}}}
\nc{\dtor}{{R-{\rm tor}}}

\nc{\Div}{{\mrm Div}} \nc{\End}{\mrm{End}} \nc{\Ext}{\mrm{Ext}}
\nc{\Fil}{\mrm{Fil}} \nc{\Frob}{\mrm{Frob}} \nc{\Gal}{\mrm{Gal}}
\nc{\GL}{\mrm{GL}} \nc{\lord}{\mrm{L-order}\xspace}
\nc{\dllord}{\mrm{DLL-order}\xspace} \nc{\rme}{\mrm{E}}
\nc{\rmt}{\mrm{T}} \nc{\Sym}{\mrm{Sym}} \nc{\Hom}{\mrm{Hom}}
\nc{\hsr}{\mrm{H}} \nc{\hpol}{\mrm{HP}} \nc{\id}{\mrm{id}}
\nc{\im}{\mrm{im}} \nc{\incl}{\mrm{incl}} \nc{\length}{\mrm{length}}
\nc{\leng}{\mrm{\ell}} \nc{\LR}{\mrm{LR}} \nc{\mchar}{\rm char}
\nc{\mzvalg}{\mathbf{MZV}} \nc{\edsalg}{\mathbf{EDS}}
\nc{\qeds}{$\QQ$-EDS\xspace} \nc{\zeds}{$\ZZ$-EDS\xspace}
\nc{\st}{stfl}
\nc{\zpeds}{$\ZZ_p$-EDS\xspace} \nc{\fpeds}{$\FF_p$-EDS\xspace}
\nc{\NC}{\mrm{NC}} \nc{\mpart}{\mrm{part}} \nc{\os}{\mrm{OS}}
\nc{\qs}{\mrm{QS}} \nc{\ql}{{\QQ_\ell}} \nc{\qp}{{\QQ_p}}
\nc{\rank}{\mrm{rank}} \nc{\rcot}{\mrm{cot}} \nc{\rdef}{\mrm{def}}
\nc{\rdiv}{{\rm div}} \nc{\rtf}{{\rm tf}} \nc{\rtor}{{\rm tor}}
\nc{\res}{\mrm{res}} \nc{\sh}{\mrm{MS}} \nc{\TL}{\mrm{TL}}
\nc{\Spec}{\mrm{Spec}} \nc{\tor}{\mrm{tor}} \nc{\Tr}{\mrm{Tr}}
\nc{\tr}{\mrm{tr}} \nc{\TEC}{\mathrm{TEC}} \nc{\TEL}{\mathrm{TEL}}
\nc{\EL}{\mathrm{EL}} \nc{\RETL}{\mathrm{RETL}}
\nc{\EETL}{\widetilde{\TL}} \nc{\word}{\rm word\xspace}
\nc{\words}{\rm words\xspace} \nc{\varab}{\phi_{\alpha,\beta}}
\nc{\lengord}{_{\mrm{leng}}} \nc{\lone}{_{\hskip -7.5pt 1}}
\nc{\pfpair}[2]{\big(\begin{array}{c}\scs{#1} \\ \scs{#2}
\end{array} \big)}

\nc{\ab}{\mathbf{Ab}} \nc{\Alg}{\mathbf{Alg}}
\nc{\Algo}{\mathbf{Alg}^0} \nc{\Bax}{\mathbf{Bax}}
\nc{\Baxo}{\mathbf{Bax}^0} \nc{\RBo}{\mathbf{RB}^0}
\nc{\BRB}{\mathbf{RB}} \nc{\Dend}{\mathbf{DD}} \nc{\bfk}{{\bf k}}
\nc{\bfone}{{\bf 1}} \nc{\base}[1]{{a_{#1}}}
\nc{\detail}{\marginpar{\bf More detail}
    \noindent{\bf Need more detail!}
    \svp}
\nc{\Diff}{\mathbf{Diff}}
\nc{\gap}{\marginpar{\bf Incomplete}\noindent{\bf Incomplete!!}
    \svp}
\nc{\FMod}{\mathbf{FMod}}
\nc{\RB}{\mathbf{RB}}
\nc{\Int}{\mathbf{Int}}
\nc{\Mon}{\mathbf{Mon}}
\nc{\remarks}{\noindent{\bf Remarks: }} \nc{\Rep}{\mathbf{Rep}}
\nc{\Rings}{\mathbf{Rings}}
\nc{\Sets}{\mathbf{Sets}}\nc{\DT}{\mathbf{DT}} \nc{\ug}{{U}}
\nc{\ssg}[1]{\overline{#1}}
\nc{\qsshab}{{{\ssha\hspace{-2pt}_{\rho}}}\,}

\nc{\nvec}[1]{[\vec{#1}]}

\nc{\BA}{{\Bbb A}}
\nc{\CC}{{\Bbb C}}
\nc{\DD}{{\Bbb D}}
\nc{\EE}{{\Bbb E}}
\nc{\FF}{{\Bbb F}}
\nc{\GG}{{\Bbb G}}
\nc{\HH}{{\Bbb H}}
\nc{\LL}{{\Bbb L}}
\nc{\NN}{{\Bbb N}}
\nc{\QQ}{{\Bbb Q}}
\nc{\RR}{{\Bbb R}}
\nc{\TT}{{\Bbb T}}
\nc{\VV}{{\Bbb V}}
\nc{\ZZ}{{\Bbb Z}}

\nc{\cala}{{\mathcal A}}
\nc{\calap}{{\mathcal A}}
\nc{\calc}{{\mathcal C}}
\nc{\cald}{{\mathcal D}}
\nc{\cale}{{\mathcal E}}
\nc{\calf}{{\mathcal F}}
\nc{\calg}{{\mathcal G}}
\nc{\calh}{{\mathcal H}}
\nc{\cali}{{\mathcal I}}
\nc{\calj}{{\mathcal J}}
\nc{\call}{{\mathcal L}}
\nc{\calm}{{\mathcal M}}
\nc{\caln}{{\mathcal N}}
\nc{\calo}{{\mathcal O}}
\nc{\calp}{{\mathcal P}}
\nc{\calr}{{\mathcal R}}
\nc{\cals}{{\mathcal S}}
\nc{\calt}{{\mathcal T}}
\nc{\calw}{{\mathcal W}}
\nc{\calx}{{\mathcal X}}
\nc{\CA}{\mathcal{A}}

\nc\indI{\mathcal{I}}

\nc{\fraka}{{\mathfrak a}}
\nc{\frakB}{{\mathfrak B}}
\nc{\frakb}{{\mathfrak b}}
\nc{\frakd}{{\mathfrak d}}
\nc{\frakF}{{\mathfrak F}}
\nc{\frakf}{{\mathfrak f}}
\nc{\frakg}{{\mathfrak g}}
\nc{\frakL}{{\mathfrak L}}
\nc{\frakm}{{\mathfrak m}}
\nc{\frakM}{{\mathfrak M}}
\nc{\frakMo}{{\mathfrak M}^0}
\nc{\frakp}{{\mathfrak p}}
\nc{\frakw}{{\mathfrak w}}
\nc{\frakx}{{\mathfrak x}}
\nc{\ox}{\overline{\frakx}}
\nc{\frakX}{{\mathfrak X}}
\nc{\fraky}{{\mathfrak y}}

\nc{\spair}[2]{\big[\begin{array}{c}\scs{#1} \\ \scs{#2} \end{array}
\big]}

\nc{\MZV}{\mrm{MZV}\xspace} \nc{\zsg}[1]{\widehat{#1}}

\nc{\li}[1]{\textcolor{blue}{Li: #1}}
\nc{\byhs}[1]{\textcolor{red}{Bingyong: #1}}
\nc{\xie}[1]{\textcolor{green}{Xie: #1}}

\nc{\rrb}[1]{[#1]} \nc{\rrrb}[1]{\{#1\}} \nc{\ideal}[1]{\langle
#1\rangle} \nc{\refl}[1]{\overline{#1}} \nc{\rrB}{{reflexive
Rota-Baxter}\xspace} \nc{\rrrB}{{radical reflexive
Rota-Baxter}\xspace} \nc{\srb}{{strict Rota-Baxter}\xspace}
\nc{\Srb}{{Strict Rota-Baxter}\xspace} \nc{\pssha}{{\star}}
\nc{\alga}{{A}}

\renewcommand\geq{\geqslant}
\renewcommand\leq{\leqslant}

\nc{\wvec}[2]{{\scriptsize{\big [ \!\!
    \begin{array}{c} #1 \\ #2 \end{array} \!\! \big ]}}}
\nc{\vr}{{\vec r}} \nc{\vs}{{\vec s}}
\nc{\lp}{\big ( } \nc{\rp}{\big ) } \nc {\msh }{\ast}
\nc\rbop{{\lc\,\,\rc}}
\nc\rbopi[1]{{\lc_{#1} \, \rc_{#1}}}

\nc{\redtext}[1]{\textcolor{red}{#1}}
\nc{\zb}[1]{\textcolor{blue}{Bin: #1}}

\newtheorem{thm}[theorem]{Theorem}
\newtheorem{ex}[theorem]{Example}
\newtheorem{rk}[theorem]{Remark}

\newcommand{\ignore}[1]{}

\newcommand{\field}[1]{{\mathbb #1}}
\newcommand{\R}{\field{R}}
\newcommand{\C}{\field{C}}

\newcommand{\Z}{\field{Z}}

\newcommand{\N}{\field{N}}
\newcommand\cutoffsum{\mathop{-\hskip -4mm\sum}\limits}

\newcommand\altcutoffsum{\mathop{-\hskip -3mm\sum}\limits}

\def\shu{\joinrel{\!\scriptstyle\amalg\hskip -3.1pt\amalg}\,}

\nc{\mop}[1]{{\mrm #1}}
\nc{\smop}[1]{{\text{#1}}}

\def \mopl#1{\mathop{\hbox{\rm #1}}\limits}

\def\cutoffint{-\hskip -10pt\int}

\def \e {{\epsilon}}
\def \restr#1{\mathstrut_{\textstyle |}\raise-6pt\hbox{$\scriptstyle #1$}}
\def \srestr#1{\mathstrut_{\scriptstyle |}\hbox to
-1.5pt{}\raise-4pt\hbox{$\scriptscriptstyle #1$}}

\def\res{{\rm res}}

\begin{document}

\title{Double shuffle relations and renormalization of multiple zeta values}
%
%========================================================================================%
\author{Li Guo}
\address{Department of Mathematics and Computer Science,
         Rutgers University,
         Newark, NJ 07102, USA}
\email{liguo@rutgers.edu}
\author{Sylvie Paycha}
\address{Laboratoire de Math\'ematiques Appliqu\'ees,
Universit\'e Blaise Pascal (Clermont II),
Complexe Universitaire des C\'ezeaux,
63177 Aubi\`ere Cedex, France}
\email{sylvie.paycha@math.univ-bpclermont.fr}
\author{Bingyong Xie}
\address{Department of Mathematics, Peking University, Beijing, 100871, P. R. China}
\email{byhsie@math.pku.edu.cn}
\author{Bin Zhang}
\address{Yangtze Center of Mathematics,
Sichuan University, Chengdu, 610064, P. R. China}
\email{zhangbin@scu.edu.cn}

%=================================================================================
\date{\today}
%==================================================================================
%\begin{document}

\begin{abstract}
In this paper we present some of the recent progresses in multiple zeta values (MZVs). We review the double shuffle relations for convergent MZVs and summarize generalizations of the sum formula and the decomposition formula of Euler for MZVs.
We then discuss how to apply methods borrowed from
renormalization in quantum field theory and from pseudodifferential calculus to partially extend the double shuffle relations to divergent MZVs.
\end{abstract}

\maketitle

%==================================================================================

%% Classification: 11M41 (MZVs), 16W30 Hopf algebra, 81T15 QFT renormalization

%\tableofcontents

\setcounter{section}{0}

%================================================================================

\section {Introduction}
The purpose of this paper is to give a survey of recent developments in multiple zeta values (MZVs). We emphasize on the double shuffle relations
which underlie the algebraic relations among the convergent MZVs, and
on renormalization methods that aim to extend the double shuffle relations to MZVs outside of the convergent range of the nested sums defining MZVs. We also provide background on double shuffle relations and renormalization, as well
 as the closely related Rota-Baxter algebras and some analytic
 tools
 in pseudodifferential calculus in view of renormalization.

\subsection{Double shuffle relations and Euler's formulas}
A {\bf multiple zeta value (MZV)} is the special value of the complex valued function
$$\zeta(s_1,\cdots,s_k)=\sum_{n_1>\cdots >n_k\geq 1} \frac{1}{n_1^{s_1}\cdots n_k^{s_k}}$$
at positive integers $s_1,\cdots,s_k$ with $s_1\geq 2$ to insure the convergence of the nested sum. MZVs are natural generalizations of the Riemann zeta values $\zeta(s)$ to multiple variables. The two variable case (double zeta values) was already studied by Euler.

MZVs in the general case were introduced 1990s with motivations from number theory~\mcite{Za}, combinatorics~\mcite{Ho0} and quantum field theory~\mcite{BK}. Since then the subject has turned into an active area of research that involves many areas of mathematics and mathematical physics~\mcite{Ca1}. Its number theoretic significance can be seen from the fact that all MZVs are periods of mixed Tate motives over $\ZZ$ and the conjecture that all periods of mixed Tate motives are rational combinations of MZVs~\mcite{Go,GM,Te}.
\smallskip

It has been discovered that the analytically defined MZVs satisfy
many algebraic relations. Further it is conjectured that these
algebraic relations all follow from the combination of two algebra
structures: the shuffle relation and the stuffle (harmonic shuffle
or quasi-shuffle) relation~\mcite{IKZ}. This remarkable conjecture
not only links the analytic study of MZVs to the algebraic study of
double shuffle relations, but also implies the more well-known
conjecture on the algebraic independence of $\zeta(2),\zeta(2k+1),
k\geq 1,$ over $\QQ$.

Many results on algebraic relations among MZVs can be regarded as generalizations of
Euler's sum formula and decomposition formula on double zeta values which preceded the general developments of multiple zeta values by over two hundred years. We summarize these results in Section~\mref{sec:euler}.
With the non-experts in mind, we first give in Section~\mref{sec:ds}
preliminary concepts and results on double shuffle relations for MZVs and the
related Rota-Baxter algebras.

\subsection{Renormalization}
Values of the Riemann zeta function at negative integers are defined by analytic continuation and possess significant number theory properties (Bernoulli numbers, Kummer congruences, $p$-adic $L$-functions, $\cdots$). Thus it would be interesting to similarly study MZVs outside of the convergent domain of the corresponding nested sums. However, most of the MZVs remain undefined even after the analytic continuation. To bring new ideas into the study, we introduce the method of renormalization from quantum field theory.

Renormalization is a process motivated by physical insight to extract finite
values from divergent Feynman integrals in quantum field theory, after
adding in a so-called counter-term.
Despite its great success in physics, this process was well-known for its lack of a solid mathematical foundation until
%The Feynman graphs appeared to be unrelated to any mathematical structure that might underlie the renormalization prescription.
%\smallskip
%
the seminal work of Connes and Kreimer~\mcite{CK,CK1,CK2,Kr}.
They obtained a Hopf algebra structure on Feynman graphs and showed that the
separation of Feynman integrals into the renormalized values and the
counter-terms comes from their algebraic Birkhoff decomposition
similar to the Birkhoff decomposition of a loop map.
\smallskip

The work of Connes and Kreimer establishes a bridge that allows an exchange of ideas between physics and mathematics. In one
direction, their work provides the renormalization of quantum field
theory with a mathematical foundation which was previously missing, opening the door to further mathematical understanding of renormalization. For example, the related Riemann-Hilbert correspondence and motivic Galois groups were studied by Connes and Marcolli~\mcite{CMa}, and motivic properties of Feynman graphs and integrals were studied by Bloch, Esnault and Kreimer~\mcite{BEK}. See~\mcite{AM,BK,Ma} for more recent studies on the motivic aspect of Feynman rules and renormalization.

In the other direction, the mathematical formulation of renormalization provided by the algebraic Birkhoff decomposition allows the method of renormalization dealing with divergent Feynman integrals in physics to be applied to divergent problems in mathematics that could not be dealt with in the past, such as the divergence in multiple zeta values~\mcite{GZ,GZ2,Zh2,MP2} and Chen symbol integrals~\mcite{MP1,MP2}.
We survey these studies on renormalization in mathematics in
Sections~\mref{sec:gz} and \mref{sec:mp} after reviewing in
Section~\mref{sec:ck} the general framework of algebraic Birkhoff
decomposition in the context of Rota-Baxter algebras. We further
  present an alternative renormalization method
  using Speer's generalized evaluators \cite{S} and  show
it leads to the same renormalized double zeta values as the algebraic Birkhoff
 decomposition method.

We hope our paper will expose this active area to a wide range of audience and
promote its further study, to gain a more thorough understanding of the double
shuffle relations for convergent MZVs and to establish a systematical
renormalization theory for the divergent MZVs. One topic that we find
of interest is to compare the various renormalization methods
presented in this paper from an abstract point of view in terms of a
renormalization group yet to be described in this context, again motivated by the study in quantum field theory. With implications back to physics in mind, we note that MZVs offer a relatively handy and tractable field of experiment for such issues when compared with the very complicated Feynman integral computations.

\medskip
\noindent
{\bf Acknowledgements:} L. Guo acknowledges the support from NSF grant DMS-0505643 and thanks JAMI at Johns Hopkins University for its hospitality. S. Paycha is grateful to D. Manchon for his comments on a preliminary version of part of this paper. B. Zhang acknowledges the support from NSFC grant 10631050.

\section{Double shuffle relations for convergent multiple zeta values}
\mlabel{sec:ds}
All rings and algebras in this paper are assumed to be unitary unless otherwise specified. Let $\bfk$ be a commutative ring whose identity is denoted by $1$.

\subsection{Rota-Baxter algebras}
Let $\lambda\in \bfk$ be fixed. A unitary (resp. nonunitary) {\bf Rota--Baxter
$\bfk$-algebra {(RBA)} of weight $\lambda$} is a pair $(R,P)$
in which $R$ is a unitary (resp. nonunitary) $\bfk$-algebra and
$P: R \to R$ is a $\bfk$-linear map such that
\begin{equation}
 P(x)P(y) = P(xP(y))+P(P(x)y)+ \lambda P(xy),\ \forall x,\ y\in R.
\mlabel{eq:Ba}
\end{equation}
In some references such as~\mcite{MP2}, the notation $\theta=-\lambda$ is used.

We will mainly consider the following Rota-Baxter operators in this paper. See~\mcite{EGsu,Gusu,Ro} for other examples.

\begin{exam} {\bf (The integration operator)}
{\rm
Define the integration operator
\begin{equation}
 I(f)(x)=\int_0^x f(t)dt
 \mlabel{eq:int}
 \end{equation}
on the algebra $C[0,\infty)$ of continuous functions $f(x)$ on $[0,\infty)$.
Then it follows from the integration by parts formula that $I$ is a Rota-Baxter operator of weight 0~\mcite{Ba}.
}
\mlabel{ex:int}
\end{exam}

\begin{exam} {\bf (The summation operator)}
{\rm
Consider the summation operator~\mcite{Zud}
$$P(f)(x):= \sum_{n\geq 1} f(x+n).$$
Under certain convergency conditions, such as $f(x) = O(x^{-2})$
and $g(x)=O(x^{-2})$,
$P(f)(x)$ and $P(g)(x)$ define absolutely convergent series and we have
\begin{eqnarray}
\lefteqn{P(f)(x) P(g)(x) = \sum_{m\geq 1} f(x+m) \sum_{n\geq 1} g(x+n)}\notag \\
&=& \sum_{n>m\geq 1} f(x+m) g(x+n) + \sum_{m>n\geq 1} f(x+m)g(x+n)
+ \sum_{m\geq 1} f(x+m)g(x+m) \mlabel{eq:sum1}\\
&=& P(f P(g))(x) + P(gP(f))(x) + P(fg)(x).\notag
\end{eqnarray}
Thus $P$ is a Rota-Baxter operator of weight 1.
}
\mlabel{ex:sum}
\end{exam}
%However, the sum operator and its iteration are not defined on some other functions. So we can only expect that subsets of a Rota--Baxter algebra can be applied to the MZV study.

\begin{exam} {\bf (The partial sum operator)}
{\rm The operator $P$ defined on sequences $\sigma:\N\to\C$ by:
\begin{equation}\mlabel{operatorP}
P(\sigma)(n)=\sum_{k=0}^{n} \sigma(k)
\end{equation}
satisfies the Rota-Baxter relation with weight $-1$. Similarly, the
operator $Q=P-Id$ which acts  on sequences $\sigma:\N\to\C$ by:
\begin{equation}\mlabel{operatorQ}
Q(\sigma)(n)=\sum_{k=0}^{n-1} \sigma(k)
\end{equation}
satisfies the Rota-Baxter relation with weight $1$.
}
\mlabel{ex:ps}
\end{exam}
\begin{exam}\label{ex:Pi} {\bf (Laurent series)}
{\rm Let $A=\bfk[\vep^{-1},\vep]]$ be the algebra of Laurent series. Define $\Pi:A\to A$ by
$$\Pi\big(\sum_{n} a_n \vep^n \big)=\sum_{n<0} a_n \vep^n.$$
Then $\Pi$ is a Rota-Baxter operator of weight $-1$.
}
\mlabel{ex:lau}
\end{exam}

\subsection{Shuffles, quasi-shuffles and mixable shuffles} \mlabel{ss:msh}
We briefly recall the construction of shuffle, stuffle and quasi-shuffle products in the framework of mixable shuffle algebras~\mcite{GK1,GK2}.

Let $\bfk$ be a commutative ring.
Let $A$ be a commutative $\bfk$-algebra {\em that is not necessarily unitary}. For a given $\lambda\in \bfk$, the {\bf mixable shuffle algebra of weight $\lambda$ generated by $A$} (with coefficients in $\bfk$) is $\sh(A)=\sh_{\bfk,\lambda}(A)$ whose underlying $\bfk$-module is that of the tensor algebra
\begin{equation}
T(A)= \bigoplus_{k\ge 0}
    A^{\otimes k}
= \bfk \oplus A\oplus A^{\otimes 2}\oplus \cdots
\mlabel{eq:mshde}
\end{equation}
equipped with the {\bf
mixable shuffle product $\shprl$ of weight $\lambda$} defined as
follows.

For pure tensors
$\fraka=a_1\ot \ldots \ot a_m\in A^{\ot m}$ and $\frakb=b_1\ot
\ldots \ot b_n\in A^{\ot n}$, a {\bf shuffle} of $\fraka$ and
$\frakb$ is a tensor list of $a_i$ and $b_j$ without change the
natural orders of the $a_i$s and the $b_j$s.
More precisely, for $\sigma\in \Sigma_{k,\ell}:=\{\tau\in S_{k+\ell}\ |\ \tau^{-1}(1)<\cdots< \tau^{-1}(k), \tau^{-1}(k+1)<\cdots<\tau^{-1}(k+\ell)\},$ the shuffle of $\fraka$ and $\frakb$ by $\sigma$ is
$$
\fraka \ssha_{\sigma} \frakb:= c_{\sigma(1)}\ot \cdots \ot c_{\sigma(k+\ell)}, \qquad \text{where } c_i=\left \{\begin{array}{ll} a_i, & 1\leq i\leq k,\\ b_{i-k}, & k+1\leq i\leq k+\ell \end{array} \right .
$$
The {\bf shuffle product} of $\fraka$ and $\frakb$ is
$$
\fraka \ssha \frakb:= \sum_{\sigma\in \Sigma_{k,\ell}} \fraka \ssha_{\sigma} \frakb.
$$

More generally, for a
fixed $\lambda\in \bfk$, a {\bf mixable shuffle} (of weight
$\lambda$) of $\fraka$ and $\frakb$ is a shuffle of $\fraka$ and
$\frakb$ in which some (or {\it none}) of the pairs $a_i\ot b_j$ are
merged into $\lambda\, a_i b_j$. Then the {\bf mixable shuffle product of weight $\lambda$} is defined by
\begin{equation}
\fraka \shprl \frakb= \sum {\rm\ mixable\ shuffles\ of\ }\fraka {\rm\ and\ } \frakb
\mlabel{eq:msh}
\end{equation}
where the subscript $\lambda$ is often suppressed when there is no danger of confusion.
For example,
$$
a_1 \shprl (b_1\ot b_2):
= \underbrace{a_1\ot b_1\ot b_2 + b_1\ot a_1\ot b_2
+ b_1\ot b_2\ot a_1}_{\rm shuffles}
+ \underbrace{\lambda (a_1 b_1)\ot b_2 + \lambda b_1\ot (a_1 b_2)}_{\rm merged\ shuffles}.
$$
With $\bfone\in \bfk$ as the unit, this product makes $T(A)$ into a commutative $\bfk$-algebra that we denote by $\sh_{\bfk,\lambda}(A)$.
See~\mcite{GK1} for further details on the mixable shuffle product.
When $\lambda=0$, we simply have the shuffle product which is also defined when $A$ is only a $\bfk$-module, treated as an algebra with the zero multiplication.

We have the following relation between mixable shuffle product and free commutative Rota-Baxter algebras.
A Rota-Baxter algebra homomorphism $f:(R,P)\to (R',P')$ between Rota-Baxter $\bfk$-algebras $(R,P)$ and $(R',P')$ is a $\bfk$-algebra homomorphism $f:R\to R'$ such that $f\circ P = P'\circ f$.

\begin{theorem} (\mcite{GK1}) The tensor product algebra $\sha(A):=\sha_{\bfk,\lambda}(A)= A\ot \sh_{\bfk,\lambda}(A)$, with the linear operator $P_A:\sha(A)\to \sha(A)$ sending $\fraka \to 1\ot \fraka$, is the free commutative Rota-Baxter algebra generated by A in the following sense. Let $j_A:A\rar \sha (A)$ be the canonical inclusion map.
Then for any
Rota-Baxter $\bfk$-algebra $(R,P)$ and any $\bfk$-algebra
homomorphism $\varphi:A\rar R$, there exists a unique Rota-Baxter
$\bfk$-algebra homomorphism $\tilde{\varphi}:(\sha (A),P_A)\rar
(R,P)$ such that $\varphi = \tilde{\varphi} \circ j_A$ as
$\bfk$-algebra homomorphisms.
\mlabel{thm:freea}
\end{theorem}

The product $\shprl$ can also be defined by the following
recursion~\mcite{EGsh,GZ2,MP2} which provides the connection between mixable shuffle algebras and quasi-shuffle algebras of Hoffman~\mcite{Ho2}.
First we define the multiplication by $A^{\ot 0}=\bfk$ to be the scalar product. In particular, $\bfone$ is the identity.
For any $m,n\geq 1$ and
$\fraka:=a_1\ot\cdots \ot a_m\in A^{\ot m}$, $\frakb:=b_1\ot \cdots\ot b_n\in A^{\ot n}$, define
$a \shprl b$ by induction on the sum $m+n\geq 2$. When $m+n=2$, we have
$a=a_1$ and $b=b_1$. We define
$$
 a\shprl b = a_1\ot b_1 + b_1\ot a_1 + \lambda a_1b_1.
$$
Assume that $\fraka\shprl \frakb$ has been defined for $m+n\geq k\geq 2$ and consider $\fraka$ and $\frakb$
with $m+n=k+1$. Then $m+n\geq 3$ and so at least one of $m$ and $n$ is greater than 1.
We define
\begin{equation}
  \fraka \shprl \frakb =\left\{ \begin{array}{l}
  a_1\ot  b_1\ot  \cdots \ot b_n  + b_1\ot \big(a_1\shprl (b_2\ot \cdots\ot b_n)\big) \\
    \qquad \qquad + \lambda(a_1b_1)\ot  b_2\ot \cdots\ot  b_n, {\rm\ when\ } m=1, n\geq 2, \\
    a_1 \ot\big ((a_2\ot \cdots\ot  a_m)\shprl b_1 \big) + b_1\ot a_1\ot \cdots\ot a_m \\
 \qquad \qquad + \lambda(a_1b_1) \ot  a_2\ot \cdots\ot  a_m,     {\rm\ when\ } m\geq 2, n=1, \\
     a_1\ot  \big ((a_2\ot \cdots\ot a_m)\shprl
(b_1\ot \cdots\ot  b_n)\big ) + b_1\ot  \big ((a_1\ot  \cdots \ot a_m)\shprl (b_2 \ot \cdots \ot  b_n)\big) \notag \\
 \qquad \qquad  + \lambda(a_1 b_1)  \big ( (a_2\ot \cdots\ot a_m) \shprl
     (b_2\ot \cdots\ot  b_n)\big ),
     {\rm\ when\ } m, n\geq 2.
     \end{array} \right .
\mlabel{eq:quasi}
\end{equation}
Here the products by $\shprl$ on the right hand side of the equation are well-defined
by the induction hypothesis.

Let $S$ be a semigroup and let $\bfk\,S=\sum_{s\in S} \bfk\,s$ be the semigroup nonunitary $\bfk$-algebra. A canonical $\bfk$-basis of $(\bfk\,S)^{\ot k}, k\geq 0$, is the set $S^{\ot k}:=\{s_1\ot \cdots \ot s_k\ |\  s_i\in S, 1\leq i\leq k\}$.
Let $S$ be a graded semigroup $S=\coprod_{i\geq 0} S_i$, $S_iS_j\subseteq S_{i+j}$ such that $|S_i|<\infty$, $i\geq 0$.
Then the mixable shuffle product $\shpr_1$ of weight $1$ is identified with the {\bf quasi-shuffle product} $\ast$ defined by Hoffman~\mcite{Ho2,EGsh,GZ2}.
\begin{notation}
{\rm To simplify the notation and to be consistent with the usual notations in the literature on multiple zeta values, we will identify $s_1\ot \cdots \ot s_k$ with the concatenation $s_1\cdots s_k$ unless there is a danger of confusion.
We also denote the weight $1$ mixable shuffle product $\shpr_1$ by $\ast$ and denote the corresponding mixable algebra $\sh_{\bfk,1}(A)$ by $\calh_A^\ast$. Similarly, when $A$ is taken to be a $\bfk$-module, we denote the weight zero mixable shuffle algebra $\sh_{\bfk,0}(A)$ by $\calh_A^\shf$.
}
\mlabel{no:mscase}
\end{notation}

Yet another interpretation of the mixable shuffle or quasi-shuffle product can be given in terms of order preserving maps that are called {\bf stuffle} in the study of MZVs but could be traced back to Cartier's work~\mcite{Ca} on free commutative Rota-Baxter algebras.

For positive integers $k$ and $\ell$, denote
$[k]=\{1,\cdots,k\}$ and $[k+1,k+\ell]=\{k+1,\cdots,k+\ell\}.$
Define
\begin{equation}
\indI_{k,\ell}=\left \{(\varphi,\psi)\ \Big|\ \begin{array}{l}
\varphi: [k]\to [k+\ell], \psi: [\ell]\to [k+\ell]
\text{ are order preserving } \\
\text{ injective maps and } \im(\varphi)\sqcup\im (\psi)=[k+\ell] \end{array} \right
\} \mlabel{eq:ind}
\end{equation}
Let $\fraka\in A^{\ot k}$, $\frakb\in A^{\ot \ell}$ and
$(\varphi,\psi)\in\indI_{k,\ell}$. We define
$\fraka\ssha_{(\varphi,\psi)}\frakb$ to be the tensor whose $i$-th
factor is
\begin{equation}
(\fraka\ssha_{(\varphi,\psi)} \frakb)_i =\left\{\begin{array}{ll}
a_j & \text{if } i=\varphi(j) \\ b_j & \text{if }
i=\psi(j)\end{array}\right.
= a_{\varphi^{-1}(i)}b_{\psi^{-1}(i)}, \quad 1\leq i\leq k+\ell,
\mlabel{eq:mulind}
\end{equation}
with the convention that $a_\emptyset=b_\emptyset=1.$
Then we have
\begin{equation}
\fraka \ssha \frakb = \sum_{(\varphi,\psi)\in \indI_{k,\ell}} \fraka \ssha_{(\varphi,\psi)} \frakb
\mlabel{eq:mapsh}.
\end{equation}

More generally, for $0\leq r \leq \min(k,\ell)$, define
$$
\indI_{k,\ell,r}=\left \{(\varphi,\psi)\ \Big|\ \begin{array}{l}
\varphi: [k]\to [k+\ell-r], \psi: [\ell]\to [k+\ell-r]
\text{ are order preserving } \\
\text{ injective maps and } \im(\varphi)\cup\im (\psi)=[k+\ell-r] \end{array} \right
\}
$$
Clearly, $\indI_{k,\ell,0}=\indI_{k,\ell}$.
Let $\fraka\in A^{\ot k}$, $\frakb\in A^{\ot \ell}$ and
$(\varphi,\psi)\in\indI_{k,\ell,r}$. We define
$\fraka\ssha_{(\varphi,\psi)}\frakb$ to be the tensor whose $i$-th factor is
$$
(\fraka\ssha_{(\varphi,\psi)} \frakb)_i =\left\{\begin{array}{ll}
a_j & \text{if } i=\varphi(j), i\notin \im\, \psi \\
b_j & \text{if } i=\psi(j), i\notin \im\, \varphi \\ a_j b_{j^\prime} & \text{if } i=\varphi(j), i=\psi(j^\prime) \end{array}\right\}
= a_{\varphi^{-1}(i)}b_{\psi^{-1}(i)}, \quad 1\leq i\leq k+\ell-r,
$$
with the convention that $a_\emptyset=b_\emptyset=1.$
Then we have~\mcite{GK1,GZ}
\begin{equation}
\fraka \shpr_\lambda \frakb = \sum_{r=0}^{\min(k,\ell)} \lambda^{r} \Big (\sum_{(\varphi,\psi)\in \indI_{k,\ell,r}} \fraka \ssha_{(\varphi,\psi)} \frakb\Big).
\mlabel{eq:mapshr}
\end{equation}
In particular,
$$
\fraka \qspr \frakb = \sum_{r=0}^{\min(k,\ell)} \Big (\sum_{(\varphi,\psi)\in \indI_{k,\ell,r}} \fraka \ssha_{(\varphi,\psi)} \frakb\Big)
=\sum_{(\varphi,\psi)\in \bar{\indI}_{k,\ell}} \fraka \ssha_{(\varphi,\psi)} \frakb
$$
where $\bar{\indI}_{k,\ell}=\cup_{r=0}^{\min(k,\ell)} \indI_{k,\ell,r}.$

Equivalently, let $\smop{\st}(k,\ell,r)$ denote the set of surjective maps from $[k+\ell]$ to $[k+\ell-r]$ that preserve the natural orders of $[k]$ and $\{k+1,\cdots,k+\ell\}$. Let
$$
\smop{\st}(k,\ell)=\bigcup_{r=0}^{\min(k,\ell)} \smop{\st}(k,\ell,r).
$$
Then
\begin{equation}
(a_1\ot\cdots\ot a_k) \qspr (a_{k+1}\ot\cdots\ot a_{k+\ell})=
\sum_{\pi\in \smop{\st}(k,\ell)} c_1^\pi\ot \cdots \ot c_{k+\ell-r}^\pi, \quad c_i^\pi = \prod_{j\in \pi^{-1}(i)} a_j.
\mlabel{eq:stuffle}
\end{equation}

A {\bf connected filtered Hopf algebra} is a Hopf
algebra $(H, \Delta)$ with $\bfk$-submodules $H^{(n)},\ n\geq 0$ of $H$ such
that
%\vspace{-.7cm}
\begin{eqnarray*}
&H^{(n)}\subseteq H^{(n+1)}, \quad
\cup_{n\geq 0} H^{(n)} = H, \quad
H^{(p)} H^{(q)}\subseteq H^{(p+q)}, &\\
&
\Delta(H^{(n)}) \subseteq \sum_{p+q=n} H^{(p)}\otimes H^{(q)}, \quad
H^{(0)}=\bfk {\rm \ (connectedness)}.&
\end{eqnarray*}

On the algebra $\sh_{\bfk,\lambda}(A)$ further define
\begin{eqnarray}
\Delta(a_1\ot \cdots \ot a_n)&=&1\bigotimes (a_1 \ot \cdots\ot a_n) + a_1 \bigotimes (a_2\ot \cdots \ot a_n) \notag \\
&&
+\cdots + (a_1\ot \cdots a_{n-1})\bigotimes a_n + (a_1\ot \cdots \ot a_n)\ot 1.
\mlabel{eq:qshco}
\end{eqnarray}
Then $\Delta$ extends by linearity to a linear map $\sh_{\bfk,\lambda}(A)\to \sh_{\bfk,\lambda}(A)\bigotimes \sh_{\bfk,\lambda}(A).$

\begin{theorem}(\mcite{GZ2,Ho2,MP2})
The triple $(\sh_{\bfk,\lambda}(A),\shprl,\Delta)$,
together with the unit $\un:\bfk \hookrightarrow \sh_{\bfk,\lambda}(A)$ and the counit $\vep:\sh_{\bfk,\lambda}(A) \to \bfk$ projecting onto the direct summand $\bfk\subseteq \sh_{\bfk,\lambda}(A)$, equips $\sh_{\bfk,\lambda}(A)$ with the structure of a connected filtered Hopf algebra with the filtration $\sh(A)^{(n)}:=\sum_{i\leq n} A^{\ot i}.$
\mlabel{thm:hopf}
\end{theorem}

%For a general semigroup $S$, the mixable shuffle algebra $\sh_{\bfk,1}(S)$ of weight $1$ coincides with the {\bf generalized overlapping shuffle algebra} on $S$~\mcite{Ha2}.

We also have the following easy extension of Hoffman's isomorphism between the shuffle Hopf algebra and the quasi-shuffle Hopf algebra (see also \mcite{EGsh}). Recall the notation $\calh_A^\ast=\sh_{\bfk,1}(A)$ and $\calh_A^\shf=\sh_{\bfk,0}(A)$.

\begin{theorem}\mlabel{th:hoffman}{\rm (\mcite{Ho2,MP2})} Let $\bfk$ be a $\QQ$-algebra. There is an isomorphism of Hopf
algebras~:
\begin{equation}
\exp:\calh_A^\shf \mathop{\longrightarrow}\limits^{\sim} \calh_A^\ast.
\mlabel{eq:his}
\end{equation}
\end{theorem}
Hoffman's isomorphism (\mref{eq:his}) is built explicitly as follows.
Let ${\mathcal P}(n)$ be the set of compositions of the integer $n$, i.e. the set of
sequences  $I=(i_1,\ldots,i_k)$ of positive integers such that
$i_1+\cdots +i_k=n$. For any  $u=v_1\otimes\cdots\otimes v_n\in T({A})$
and any composition  $I=(i_1,\ldots,i_k)$ of $n$ we set:
$$I[u]:=(v_1\cdot\cdots\cdot v_{i_1})\otimes
(v_{i_1+1}\cdot\cdots\cdot v_{i_1+i_2})
\otimes\cdots\otimes
(v_{i_1+\cdots+i_{k-1}+1}\cdot\cdots\cdot v_n).$$
Then the isomorphism $\exp$ is defined by
$$\exp u=\sum_{I=(i_1,\cdots ,i_k)\in{\mathcal P}(n)}\frac{1}{i_1!\cdots
  i_k!}I[u].$$
Moreover (\mcite{Ho2}, Lemma 2.4), the inverse $\log$ of $\exp$ is given by~:
$$\log u=\sum_{I=(i_1,\cdots ,i_k)\in{\mathcal P}(n)}\frac{(-1)^{n-k}}{i_1\cdots
  i_k}I[u].$$

\subsection{Double shuffle of MZVs and related conjectures}
\mlabel{ss:prel}
A multiple zeta value (MZV) is defined to be
\begin{equation}
 \zeta(s_1,\cdots,s_k):= \sum_{n_1>\cdots>n_k\geq 1} \frac{1}{n_1^{s_1}\cdots n_k^{s_k}}
 \mlabel{eq:sumrep}
\end{equation}
where $s_i\geq 1$ and $s_1>1$ are integers.
As is well-known, an MZV has an integral representation due to Kontsevich~\mcite{LM}
\begin{equation}
\zeta(s_1,\cdots,s_k)=
 \int_0^1 \int_0^{t_1}\cdots \int_0^{t_{|\vec{s}|-1}} \frac{dt_1}{f_1(t_1)}
 \cdots \frac{dt_{|\vec{s}|}}{f_{|\vec{s}|}(t_{|\vec{s}|})}
 \mlabel{eq:intrep}
\end{equation}
Here $|\vec{s}|=s_1+\cdots +s_k$ and
$$f_j(t)=\left\{\begin{array}{ll} 1-t_j, & j= s_1,s_1+s_2,\cdots, s_1+\cdots +s_k,\\
t_j, & \text{otherwise}. \end{array} \right .
$$
The MZVs spanned the following $\QQ$-subspace of $\RR$
$$
\mzvalg: = \QQ \{ \zeta(s_1,\cdots,s_k)\ |\ s_i\geq 1, s_1\geq 2\} \subseteq \RR.
$$

Since the summation operator in Eq.~(\mref{eq:sumrep}) and the integral
operator in Eq.~(\mref{eq:int}) are both Rota-Baxter operators (of weight 1
and 0 respectively) by Example~\mref{ex:sum} and Example~\mref{ex:int}, it can
be expected that the multiplication of two MZVs follows the multiplication
rule in a free Rota-Baxter algebra and thus in a mixable shuffle algebra. This
viewpoint naturally leads  to the following double shuffle relations of MZVs.

For the sum representation of MZVs in Eq.~(\mref{eq:sumrep}), consider the
semigroup
$$ Z:=\{z_s\ |\ s\in \ZZ_{\geq 1},\quad z_s \cdot z_t=z_{s+t}, s,t \geq
1.\}$$
With the convention in Notation~\mref{no:mscase},
  we denote the quasi-shuffle algebra $\calh\qsh:=\calh\qsh_{\hspace{-.2cm}\QQ\, Z}$ which contains the subalgebra
$$
\calh\qsh\lzero:=\QQ \oplus \Big (\bigoplus_{s_1>1} \QQ z_{s_1}
\cdots z_{s_k} \Big ).
$$
Then the multiplication rule of MZVs according to their summation
representation follows from the fact that the linear map
\begin{equation}
\zeta\qsh: \calh\qsh\lzero \to \mzvalg, \quad z_{s_1,  \cdots, s_k} \mapsto \zeta(s_1,\cdots,s_k)
\mlabel{eq:mzvst}
\end{equation}
is an algebra homomorphism~\mcite{Ho1,IKZ}.

For the integral representation of MZVs in Eq.~(\mref{eq:intrep}), consider
the set $X=\{x_0,x_1\}$. With the convention in Notation~\mref{no:mscase}, we denote the shuffle algebra $\calh\shf:=\calh\shf_{\hspace{-.2cm}\QQ\,X}$ which contains subalgebras
$$
 \calh\shf\shzero:=\QQ \oplus x_0\calh\shf x_1
\subseteq \calh\shf\lone:=\QQ \oplus \calh\shf x_1 \subseteq
\calh\shf.
$$
Then the multiplication rule of MZVs according to their integral
representations follows from the statement that the linear map
$$
 \zeta\shf: \calh\shf\shzero \to \mzvalg,
 \quad x_0^{  s_1-1}  x_1  \cdots   x_0^{s_k-1}  x_1 \mapsto \zeta(s_1,\cdots,s_k)
$$
is an algebra homomorphism~\mcite{Ho1,IKZ}.

There is a natural bijection of $\QQ$-vector spaces (but {\em not} algebras)
$$
\shqs: \calh\shf\lone \to \calh\qsh, \quad 1\leftrightarrow 1,\
x_0^{s_1-1}  x_1  \cdots   x_0^{s_k-1}  x_1
 \leftrightarrow z_{s_1, \cdots, s_k}.
$$
that restricts to a bijection of vector spaces
$\shqs: \calh\shf\shzero \to \calh\qsh\lzero.
$
Then the fact that MZVs can be multiplied in two ways is reflected
by the commutative diagram
$$
\xymatrix{ \calh\shf\shzero \ar^{\shqs}[rr] \ar_{\zeta\shf}[rd] &&
    \calh\qsh\lzero \ar^{\zeta\qsh}[ld]\\
& \mzvalg & }
$$

Through $\shqs$, the shuffle product $\ssha$ on $\calh\shf\lone$ and
$\calh\shf\shzero\ $ transports to a product $\qssha$ on $\calh\qsh$
and $\calh\qsh\lzero$. That is, for $w_1,w_2\in \calh\qsh\lzero$,
define
\begin{equation}
w_1 \qssha w_2: = \shqs( \shqs^{-1}(w_1)\ssha \shqs^{-1}(w_2)).
\mlabel{eq:shtrans}
\end{equation}
Then the {\bf double shuffle relation} is simply the set
$$
\{ w_1\qssha w_2 - w_1 \ast w_2\ |\ w_1,w_2\in\calh\qsh\lzero\}$$
and the {\bf extended double shuffle relation}~\mcite{IKZ,Ra,Zud} is the set
\begin{equation}
\{ w_1\qssha w_2 - w_1 \ast w_2,\ z_1 \qssha w_2 - z_1 \ast w_2\ |\
w_1,w_2\in\calh\qsh\lzero\}.
\mlabel{eq:eds}
\end{equation}
%The following result is well-known.
\begin{theorem} {\bf (\mcite{Ho1,IKZ,Ra})}
Let $I_\edsalg$ be the ideal of $\calh\qsh\lzero$ generated by the extended double shuffle relation in Eq.~(\mref{eq:eds}). Then $I_\edsalg$ is in the kernel of $\zeta\qsh$.
\mlabel{thm:eds}
\end{theorem}
It is conjectured that $I_\edsalg$ is in fact the kernel of $\zeta\qsh$.
A consequence of this conjecture is the irrationality of $\zeta(2n+1), n\geq 1$.

\section{Generalizations of Euler's formulas}
\mlabel{sec:euler}
We begin with stating Euler's sum and decomposition formulas in Section~\mref{ss:euler}. Generalizations of Euler's sum formula are presented in Section~\mref{ss:sum} and generalizations of Euler's decomposition formula are presented in Section~\mref{ss:dec}.

\subsection{Euler's sum and decomposition formulas}
\mlabel{ss:euler}
Over two hundred years before the general study of multiple zeta values was started in the 1990s, Goldback and Euler had already considered the two variable case, the double zeta values~\mcite{Eu,Sa}
$$ \zeta(s_1,s_2):= \sum_{n_1>n_2\geq 1} \frac{1}{n_1^{s_1} n_2^{s_2}}.
$$
Among Euler's major discoveries on double zeta values are his {\bf
sum formula} $$ \sum_{i=2}^{n-1} \zeta(i,n-i)=\zeta(n)$$ expressing
one Riemann zeta values as a sum of double zeta values and the
{\bf decomposition formula}
\begin{equation}
\zeta(r) \zeta(s)
 = \sum_{k=0}^{s-1} \binc{r+k-1}{k} \zeta(r+k,s-k)
+ \sum_{k=0}^{r-1} \binc{s+k-1}{k} \zeta(s+k,r-k), \quad r,s\geq 2,
\mlabel{eq:euler}
\end{equation}
expressing the product of two Riemann zeta values as a sum of double zeta values.

A major aspect of the study of MZVs is to find algebraic and linear relations among MZVs, such as Euler's formulas. Indeed a large part of this study can be viewed as  generalizations of Euler's formulas.

\subsection{Generalizations of Euler's sum formula}
\mlabel{ss:sum}
Soon after MZVs were introduced,
Euler's sum formula was generalized to MZVs~\mcite{Ho0,Gr,Za2} as the well-known sum formula, followed by quite a few other generalizations that we will next summarize.

\subsubsection{Sum formula} \label{ss:sumform}
The first generalization of Euler's sum formula is the sum formula
conjectured in~\mcite{Ho0}.
Let
\begin{equation}
I(n,k)=\{(s_1,\cdots, s_k)\ |\ s_1+\cdots+ s_k=n, s_i\geq 1,
s_1\geq 2\}.
\mlabel{eq:zrange}
\end{equation}
For $\vec{s}=(s_1,\cdots,s_k)\in I(n,k)$, define the {\bf multiple
zeta star value (or non-strict MZV)}
\begin{equation}
\zeta^\star(s_1,\cdots, s_k)=\sum_{n_1\geq \cdots\geq n_k\geq 1} \frac{1}{n_1^{s_1}\cdots n_k^{s_k}}.
\mlabel{eq:mszv}
\end{equation}
Note the subtle different between the notations $\zeta^\ast$ in Eq.~(\mref{eq:mzvst}) and $\zeta^\star$ in Eq.~(\mref{eq:mszv}).

\begin{theorem}
$(${\bf Sum formula}$)$ For positive integers $k <n$ we have
\begin{equation}\label{eq:sumform}
\sum\limits_{\vec{s}\in I(n,k)}\zeta(\vec{s})=\zeta(n),\hskip 10pt
\sum\limits_{\mathbf{k}\in
I(n,k)}\zeta^\star(\vec{s})=\binc{n-1}{k-1}\zeta(n).\end{equation} \mlabel{thm:sum}
\end{theorem}

The case of $k=3$ was proved by M. Hoffman and C. Moen \mcite{HM} and the general
case was proved by Zagier~\mcite{Za2} with another proof given by
Granville~\mcite{Gr}. Later S. Kanemitsu, Y. Tanigawa, M. Yoshimoto
\cite{KTY} gave a proof for the case of $k=2$ using Mellin
transformation.

J.-I. Okuda and K. Ueno \cite{OU} gave the following version of the sum
formula
$$\sum_{k=r}^n\binc{k-1}{r-1}\Big(\sum\limits_{\vec{s}\in I(n,k)}\zeta(\vec{s})\Big)=\binc{n-1}{r}\zeta(n)$$
for $n>r\geq 1$ from which they deduced the sum formula
Eq.~(\ref{eq:sumform}).

\subsubsection{Ohno's generalized duality theorem}

Another formula conjectured in~\mcite{Ho0} is the {\bf  duality formula}.
To state the duality formula, we need an involution $\tau$ on the
set of finite sequences of positive integers whose first element is
greater than $1$. If $$\vec{s}=(1+b_1,\underbrace{1,\cdots,
1}_{a_1-1},\cdots, 1+b_k, \underbrace{1,\cdots, 1}_{a_k-1}),$$ then
$$\tau(\vec{s})=(1+a_k,\underbrace{1,\cdots,
1}_{b_k-1},\cdots, 1+a_1, \underbrace{1,\cdots, 1}_{b_1-1}).$$

\begin{theorem} {\bf (Duality  formula)}
$$\zeta(\vec{s})=\zeta(\tau(\vec{s})).$$
\mlabel{thm:dual}
\end{theorem}
This formula is an immediate consequence of the integral representation in Eq.~(\mref{eq:intrep}).

Y. Ohno \mcite{Oh} provided a generalization of both the sum formula
and the duality formula.
\begin{theorem} $(${\bf Generalized Duality Formula~\mcite{Oh}}$)$
For any index set $\vec{s}=(s_1,\cdots,
s_k)$ with $s_1\geq 2,s_2\geq 1, \cdots, s_k\geq 1 $, and a
nonnegative integer $\ell$, set
$$Z(s_1,\cdots, s_k; \ell)=\sum_{\footnotesize \begin{array}{c} c_1+\cdots+c_k=\ell \\ c_i\geq 0 \end{array} }
\zeta(s_1+c_1,\cdots, s_k+c_k).$$
Then
$$Z(\vec{s}; \ell)=Z(\tau(\vec{s}); \ell).$$
\end{theorem}
When $\ell=0$,
this is just the duality formula. When $\vec{s}=(k+1)$ and
$\ell=n-k-1$, this becomes the sum formula.

\subsubsection{Sum formulas with further conditions on the variables}

M. Hoffman and Y. Ohno~\mcite{HO} gave a cyclic generalization of the sum formula.
\begin{theorem} {\bf (Cyclic sum formula)}
For any positive integers $s_1,\cdots,s_k$ with some $s_i\geq 2$,
$$
\sum_{i=1}^k \zeta(s_i+1,s_{i+1},\cdots,s_k,s_1,\cdots,s_{i-1}) = \sum_{\{i\, |\, s_i\geq 2\}} \sum_{j=0}^{s_i-2} \zeta(s_i-j,s_{i+1},\cdots,s_k,s_1,\cdots,s_{i-1},j+1).$$
\end{theorem}

Y. Ohno and N. Wakabayashi \cite{OW} gave a cyclic sum formula for non-strict MZVs and
used it to prove the sum formula Eq.~(\ref{eq:sumform}).

\begin{theorem} $(${\bf Cyclic sum formula in the non-strict case}$)$
For positive integers $k<n$ and $(s_1,\cdots, s_k)\in I(n,k)$ we
have
\begin{equation}
\sum_{i=1}^k\sum_{j=0}^{s_i-2}\zeta^\star(s_i-j,s_{i+1},\cdots,s_k, s_1,\cdots, s_{i-1}, j+1
)=n\zeta(n+1),
\mlabel{eq:csum}
\end{equation}
 where the empty sums are zero.
\end{theorem}

M. Eie, W.-C. Liaw and Y. L. Ong \mcite{ELO} gave a generalization of the sum formula by allowing a more general form in the arguments in the MZVs.

\begin{theorem}
For all positive integers $n,k$
with $n>k$, and a nonnegative integer $p$,
$$\sum_{\footnotesize\begin{array}{c}s_1+\cdots+s_k=n  \\
s_1\geq 2 \end{array}}\zeta(s_1,\cdots, s_k, \{1\}^p) \ =
\sum_{\footnotesize\begin{array}{c} c_1+\cdots+c_{p+1}=n+p\\ c_1\geq
n-k+1
\end{array}}\zeta(c_1,\cdots, c_{p+1}).$$
\end{theorem}
When $p=0$, this becomes the sum formula.
\smallskip

Y. Ohno and D.  Zagier \cite{OZa} studied another kind of sum with
certain restrictive conditions. Let
$$I(n,k,r)=\{(s_1,\cdots, s_k) \ |\ s_i\in \ZZ_{\geq 1}, s_1+\cdots+ s_k=n, \#\{s_i \ |\ s_i\geq 2
\}=r\}$$ and put
$$G(n,k,r)=\sum_{\vec{s}\in I(n,k,r)}\zeta(\vec{s}).$$ They
studied the associated generating function
$$\Phi(x,y,z)=\sum_{r\geq 1, k\geq r, n\geq k+r}G(n,k,r) x^{n-k-r}y^{k-r}z^{r-1}\in \RR[x,y,z]$$
and proved the following

\begin{theorem} We have
$$\Phi(x,y,z)=\frac{1}{xy-z}\left(1-\exp\left(\sum_{n=2}^\infty\frac{\zeta(n)}{n}S_n(x,y,z)\right)\right),$$
where $S_n(x,y,z)$ are given by the identity
$$\log\left(1-\frac{xy-z}{(1-x)(1-y)}\right)=\sum_{n=2}^\infty\frac{S_n(x,y,z)}{n}$$
and the requirement that $S_n(x,y,z^2)$ is a homogeneous polynomial
of degree $n$. In particular, all of the coefficients $G(n,k,r)$
can be expressed as polynomials in $\zeta(2), \zeta(3), . . .$ with
rational coefficients.
\end{theorem}

\subsubsection{Sum formulas for $q$-MZVs}
The concept of $q$-multiple zeta values ($q$-MZVs, or multiple $q$-zeta values) was introduced as a ``quantumization" of MZVs that recovers MZVs when $q\mapsto 1$~\mcite{Br,Zh}.

For positive integers $s_1,\cdots,s_k$ with $s_1\geq 2$, define the $q$-MZV
$$\zeta_q(s_1,\cdots, s_k)=\sum_{n_1> \cdots> n_k\geq 1}
\frac{q^{n_1(s_1-1)+\cdots+ n_k(s_k-1) }}{[n_1]^{s_1}\cdots
[n_k]^{s_k}}$$
and the non-strict $q$-MZV
$$\zeta_q^\star(s_1,\cdots, s_k)=\sum_{n_1\geq \cdots\geq n_k\geq 1}
\frac{q^{n_1(s_1-1)+\cdots+ n_k(s_k-1)}}{[n_1]^{s_1}\cdots
[n_k]^{s_k}},$$ where $[n]=\frac{1-q^{n}}{1-q}$.

D.~M.~Bradley \cite{Br} proved the $q$-analogue of the sum formula for $\zeta_q$.

\begin{theorem}\label{th:q-analogy-Br}  $(${\bf $q$-analogue of the sum formula}$)$
For positive integers $0<k<n$ we have
\begin{equation} \label{eq:q-sum-form}
\sum_{\substack{s_i\geq 1, s_1\geq 2\\
s_1+\cdots+s_{k}=n}}\zeta_q(s_1,\cdots,s_k)=\zeta_q(n).
\end{equation}
\end{theorem}

Y. Ohno and J.-I. Okuda \cite{OO} gave the following $q$-analogue of
the cyclic sum formula~(\mref{eq:csum}) and then used it to prove a
$q$-analogue of the sum formula for
$\zeta_q^\star$.

\begin{theorem} $(${\bf $q$-analogue of the cyclic sum formula}$)$
For positive integers $0<k<n$ and $(s_1,\cdots, s_k)\in I(n,k)$ we
have
$$\sum_{i=1}^k\sum_{j=0}^{s_i-2}\zeta^\star_q(s_i-j,s_{i+1},\cdots,s_k, s_1,\cdots, s_{i-1}, j+1
)=\sum_{\ell=0}^{k}(n-\ell)\binc{k}{\ell}(1-q)^\ell\zeta_q(n-\ell+1),$$
where the empty sums are zero.
\end{theorem}

\begin{theorem}\label{th:q-analogy}  $(${\bf $q$-analogue of the sum formula in the non-strict case}$)$
For positive integers $0<k<n$ we have
\begin{equation} \label{eq:q-sum-forms}
\sum_{\substack{s_i\geq 1, s_1\geq 2\\
s_1+\cdots+s_{k}=n}}\zeta^\star_q(s_1,\cdots,s_k)=\frac{1}{n-1}\binc{n-1}{k-1}
\sum_{\ell=0}^{k-1}(n-1-\ell)(1-q)^\ell\zeta_q(n-\ell).
\end{equation}
\end{theorem}

% Other generalizations of the sum formula can be found in~\mcite{Br,KTY,OZa,OU}.

\subsubsection{Weighted sum formulas}
In the other direction to generalize Euler's sum formula, there is the weighted version of Euler's sum formula recently obtained by Ohno and Zudilin~\mcite{OZ}.

\begin{theorem}
$(${\bf Weighted Euler's sum formula}~\mcite{OZ}$)$ For any integer
$n\geq 2$, we have
\begin{equation}
\sum_{i=2}^{n-1} 2^i \zeta(i,n-i) = (n+1)\zeta (n).
\mlabel{eq:oz}
\end{equation}
\mlabel{thm:wes}
\end{theorem}
They applied it to study multiple zeta  star values. By the sum
formula, Eq.~(\mref{eq:oz}) is equivalent to the following equation
\begin{equation}
\sum_{i=2}^{n-1} (2^i-1) \zeta(i,n-i) = n\zeta (n). \label{eq:oz1}
\end{equation}

As a generalization of Eq.~(\ref{eq:oz1}), two of the authors proved
the following

\begin{theorem}
$(${\bf Weighted sum formula~\mcite{GX3}}$)$ For positive integers
$k\geq 2$ and $n\geq k+1$, we have
$$
\sum_{\substack{s_i\geq 1, s_1\geq 2\\
s_1+\cdots+s_{k}=n}}\hskip -10pt
\Big[2^{s_1-1}+(2^{s_1-1}-1)\Big(\big(\sum_{i=2}^{k-1}
2^{S_i-s_1-(i-1)}\big)+ 2^{S_{k-1}-s_1-(k-2)}\Big)\Big]
\zeta(s_1,\cdots, s_{k})=n\zeta(n),
$$
where $S_i=s_1+\cdots+s_i$ for $i=1,\cdots, k-1$.
\mlabel{thm:wsf}
\end{theorem}

\subsection{Generalizations of Euler's decomposition formula}
\mlabel{ss:dec}
Unlike the numerous generalizations of Euler's sum formula, no generalization of Euler's decomposition formula to MZVs,
neither proved nor conjectured,  had been given until~\mcite{GX2}
even though Euler's decomposition formula was recently revisited in connection with modular forms~\mcite{GKZ} and weighted sum formula~\mcite{OZ} on weighted sum formula of double zeta values, and was generalized to the product of two $q$-zeta values~\mcite{Br2,Zh}.
%
%This lack of progress should be due to the lack of a suitable context in which a generalization of Euler's formula makes sense, and the correct language and tools to formulate and prove such a generalization.

\subsubsection{Euler's decomposition formula and double shuffle}

The first step in generalizing Euler's decomposition formula is to place it as a special case in a suitable broader context. In~\mcite{GX2}, Euler's decomposition formula was shown to be a special case of the double shuffle relation. We give a proof of Euler's formula in this context before presenting its generalization in the next subsection.

We recall that the extended double shuffle relation is the set
$$
\{ w_1\qssha w_2 - w_1 \ast w_2,\ z_1 \qssha w_2 - z_1 \ast w_2\ |\
w_1,w_2\in\calh\qsh\lzero\}.
$$
Thus the determination of the double shuffle relation amounts to computing the two products $\ast$ and $\qssha$.

It is straightforward to compute the product $\ast$, either from its recursive definition in Eq.~(\mref{eq:quasi}) or its explicit interpretation as mixable shuffles in Eq.~(\mref{eq:msh}) and stuffles in Eq.~(\mref{eq:mapshr}) or (\mref{eq:stuffle}).
For example, to determine the double shuffle relation from
multiplying two Riemann zeta values $\zeta(r)$ and $\zeta(s)$,
$r,s\geq 2$, one uses their sum representations and easily gets the
quasi-shuffle relation
$$
 \zeta(r)\zeta(s)=\zeta(r,s)+\zeta(s,r)+\zeta(r+s).
$$

On the other hand, computing the product $\qssha$ is more involved as can already be seen from its definition in Eq.~(\mref{eq:shtrans}).
One first needs to use
their integral representations to express $\zeta(r)$ and $\zeta(s)$
as iterated integrals of dimensions $r$ and $s$, respectively. One
then uses the shuffle relation to express the
product of these two iterated integrals as a sum of $\binc{r+s}{r}$
iterated integrals of dimension $r+s$. Finally, these last iterated
integrals are translated back to MZVs and give the shuffle relation of $\zeta(r)\zeta(s)$.
As an illustrating example, consider $\zeta(100)\zeta (200)$. The quasi-shuffle relation is simply $\zeta(100)\zeta(200)=\zeta(100,200)+\zeta(200,100)+\zeta(300)$, but the shuffle relation is a large sum of $\binc{300}{100}$ shuffles of length (dimension) 300.
As we will show below, an explicit formula for this is precisely Euler's decomposition formula~(\mref{eq:euler}).

\begin{theorem}
For $r,s\geq 2$, we have
$$
 z_r \qssha z_s  = \sum_{k=0}^{s-1} \binc{r+k-1}{k} z_{r+k}z_{s-k}
+ \sum_{k=0}^{r-1} \binc{s+k-1}{k} z_{s+k}z_{r-k}.
$$
\mlabel{thm:eds1}
\end{theorem}
Via the algebra homomorphism $\zeta^\ast$ in Eq.~(\mref{eq:mzvst}) this theorem immediately
gives Euler's decomposition formula. Applying to the above example,
we have
$$ \zeta(100)\zeta(200)=\sum_{k=0}^{199} \binc{100+k-1}{k} \zeta(100+k,200-k) + \sum_{k=0}^{99} \binc{200+k-1}{k} \zeta(200+k,100-k).$$
\begin{proof}
Following the definition of $\qssha$ in Eq.~(\mref{eq:shtrans}), we have
$$ z_r \qssha z_s = \shqs( x_0^{r-1}x_1 \ssha x_0^{s-1}x_1).$$
So we just need to prove
$$x_0^{r-1}x_1 \ssha x_0^{s-1}x_1 =
\sum_{k=0}^{s-1} \binc{r+k-1}{k} x_0^{r+k-1}x_1 x_0^{s-k-1}x_1
+ \sum_{k=0}^{r-1} \binc{s+k-1}{k} x_0^{s+k-1}x_1 x_0^{r-k-1}x_1$$
since $\qssha(x_0^{r+k-1}x_1 x_0^{s-k-1}x_1)=z_{r+k}z_{s-k}$ and
$\qssha(x_0^{s+k-1}x_1 x_0^{r-k-1}x_1)=z_{s+k}z_{r-k}.$
This has a direct shuffle proof~\mcite{3BL2}. But we use the description of order preserving maps of shuffles in order to motivate the general case.

By Eq.~(\mref{eq:mapsh}), we have
$$ x_0^{r-1}x_1 \ssha x_0^{s-1}x_1 =
\sum_{(\varphi,\psi)\in \indI(r,s)} x_0^{r-1}x_1 \ssha_{(\varphi,\psi)} x_0^{s-1}x_0.$$
Since $\varphi$ and $\psi$ are order preserving, we have the disjoint union $\indI(r,s)= \indI(r,s)' \sqcup \indI(r,s)''$ where
$$\indI(r,s)'= \{ (\varphi,\psi)\in \indI(r,s)\ |\ \psi(s)=r+s\}$$
and
$$\indI(r,s)''=\{(\varphi,\psi)\in \indI(r,s)\ |\ \varphi(r)=r+s\}.$$
Again by the order preserving property, for $(\varphi,\psi)\in \indI(r,s)'$, we must have $\varphi(r)=r+k$ where $k\geq 0$. Thus for such $(\varphi,\psi)$, we have
$$x_0^{r-1}x_1 \ssha_{(\varphi,\psi)} x_0^{s-1} x_1 =
x_0^{r-1+k}x_1x_0^{s-1-k}x_1$$
since $\im\,\varphi \sqcup \im\, \psi =[r+s].$
For fixed $k\geq 0$, $\varphi(r)=r+k$ means that there are $k$ elements $i_1,\cdots,i_k$ from $[s-1]$ such that $\psi(i_j)\in [r+k-1]$ since $\psi(s)=r+s$. Thus $k\geq s-1$ and, since $\psi$ is order preserving, we have $\{i_1,\cdots,i_k\}=[k]$. Further there are $\binc{r+k-1}{k}$ such $\psi$'s since $\psi([k])$ can take any $k$ places in $[r+k-1]$ in increasing order and then $\phi([r])$ takes the rest places in increase order.
Thus
$$ \sum_{(\varphi,\psi)\in \indI(r,s)'} x_0^{r-1}x_1 \ssha_{(\varphi,\psi)} x_0^{s-1}x_1 = \sum_{k= 0}^{s-1} \binc{r+k-1}{k} x_0^{r+k-1}x_1 x_0^{s-k-1}x_1.$$
By a similar argument, we have
$$ \sum_{(\varphi,\psi)\in \indI(r,s)''} x_0^{r-1}x_1 \ssha_{(\varphi,\psi)} x_0^{s-1}x_1 = \sum_{k= 0}^{r-1} \binc{s+k-1}{k} x_0^{s+k-1}x_1 x_0^{r-k-1}x_1.$$
This completes the proof.
\end{proof}

\subsubsection{Generalizations of Euler's decomposition formula}

In a recent work~\mcite{GX2}, two of the authors generalized Euler's decomposition formula in two
directions, from the product of one variable functions to that of
multiple variables and from multiple zeta values to multiple
polylogarithms.

A {\bf multiple polylogarithm value} ~\mcite{3BL,Go,Go2} is defined
by
$$
\Li_{s_1,\cdots,s_k}(z_1,\cdots,z_k):= \sum_{n_1>\cdots>n_k\geq 1}
\frac{z_1^{n_1}\cdots z_k^{n_k}}{n_1^{s_1}\cdots n_k^{s_k}}
$$
where $|z_i|\leq 1$, $s_i\in \ZZ_{\geq 1}$, $1\leq i\leq k$, and
$(s_1,z_1)\neq (1,1)$. When $z_i=1, 1\leq i\leq k$, we obtain the
multiple zeta values $
 \zeta(s_1,\cdots,s_k)$. With the notation of~\mcite{3BL}, we have
$$
\begin{aligned}
&\Li_{s_1,\cdots,s_k}(z_1,\cdots,z_k)= \lambda\big(\begin{array}{c}
s_1,\cdots,s_k\\ b_1,\cdots,b_k \end{array} \big):=
\sum_{n_1>n_2\cdots >n_k\geq 1}\frac{\big(\frac{1}{b_1}\big)^{n_1}
\big(\frac{b_1}{b_2}\big)^{n_2}\cdots
\big(\frac{b_{k-1}}{b_k}\big)^{n_k}} {n_1^{s_1}n_2^{s_2}\cdots
n_k^{s_k}},
\end{aligned}
$$
where $(b_1,\cdots,b_k)= (z_1^{-1}, (z_1z_2)^{-1}, \cdots,
(z_1\cdots z_k)^{-1})$.

To state the result, let $k$ and $\ell$ be positive integers and let $\indI_{k,\ell}$ be as defined in Eq.~(\mref{eq:ind}).
Let $\vec{r}=(r_1,\cdots, r_k)\in\ZZ_{\geq 1}^k$,
$\vec{s}=(s_1,\cdots,s_\ell)\in \ZZ_{\geq 1}^{\ell}$ and
$\vec{t}=(t_1,\cdots, t_{k+\ell})\in \ZZ_{\geq 1}^{k+\ell}$ with
$|\vec{r}|+|\vec{s}|=|\vec{t}|$. Here $|\vec{r}|=r_1+\cdots +r_k$
and similarly for $|\vec{s}|$ and $|\vec{t}|$.
Denote $R_i=r_1+\cdots +r_i$ for $i\in [k]$, $S_i=s_1+\cdots +s_i$
for $i\in [\ell]$ and $T_i=t_1+\cdots+t_i$ for $i\in [k+\ell]$.
For $(\varphi,\psi)\in \indI_{k,\ell}$ and $i\in [k+\ell]$, define
$$
h_{(\varphi,\psi),i}=h_{(\varphi,\psi),(\vec{r},\vec{s}),i} =
       \left\{
              \begin{array}{ll} r_{j} & \text{ if } i=\varphi(j)
                                 \\
                                s_{j} & \text{ if } i=\psi(j)
              \end{array}
              \right.
              = r_{\varphi^{-1}(i)}s_{\psi^{-1}(i)},
$$
with the convention that $r_\emptyset =s_\emptyset =1.$

With these notations, we define
\begin{equation}
c_{\vec{r},\vec{s}}^{\vec{t},(\varphi,\psi)}(i) =\left\{
\begin{array}{ll}
\binc{t_i-1}{h_{(\varphi,\psi),i}-1} &
    \begin{array}{l}\text{if } i=1,  \text{if }
i-1,i \in \im(\varphi)\\ \text{or if } i-1,i \in \im(\psi),
\end{array}
\vspace{.2cm}
\\ \vspace{.2cm}
\begin{array}{l}
\binc{t_i-1} {T_i-R_{|\varphi^{-1}([i])|}-S_{|\psi^{-1}([i])|}}\\
= \binc{t_{i}-1}{\sum\limits_{j=1}^{i} t_j -\sum\limits_{j=1}^{i}
h_{(\varphi,\psi),j}}
\end{array} & \text{ otherwise}.
\end{array}
\right. \mlabel{eq:coef-re-def1}
\end{equation}

For $\vec{a}\in (S^1)^k$ and $\vec{b}\in (S^1)^\ell$, as in Eq.~(\mref{eq:mulind}), define
\begin{equation}
\vec{a}\ssha_{(\varphi,\psi)} \vec{b} =(a_{\varphi^{-1}(1)}b_{\psi^{-1}(1)},\cdots, a_{\varphi^{-1}(k+\ell)}b_{\psi^{-1}(k+\ell)}).
 \mlabel{eq:mulind2}
\end{equation}

\begin{theorem} $($\mcite{GX2}$)$ \mlabel{thm:shLi}
Let $k,\ell$ be positive integers. Let $\vec{r}\in \ZZ_{\geq 1}^k$
and $\vec{s}\in \ZZ_{\geq 1}^\ell$. Let $\vec{a}=(a_1,\cdots,a_k)\in
(S^1)^k$ and $\vec{b}=(b_1,\cdots,b_\ell)\in (S^1)^\ell$ such that
$\spair{r_1}{a_1}\neq \spair{1}{1}$ and $\spair{s_1}{b_1}\neq
\spair{1}{1}$. Then
$$
\lambda\big(\begin{array}{c}\vec{r}\\ \vec{a}\end{array} \big)\,
\lambda\big(\begin{array}{c}\vec{s}\\ \vec{b}\end{array} \big) =
\sum_{\vec{t}\in \ZZ_{\geq 1}^{k+\ell}, |\vec{t}|=
  |\vec{r}|+|\vec{s}|}
  \sum _{(\varphi,\psi)\in \indI_{k,\ell}}
  \Big(\prod_{i=1}^{k+\ell}c_{\vec{r},\vec{s}}^{\vec{t},
(\varphi,\psi)}(i)\Big)
    \lambda\big(\begin{array}{c} \vec{t}\\ \vec{a} \ssha_{(\varphi,\psi)} \vec{b}\end{array} \big).
$$
where $c_{\vec{r},\vec{s}}^{\vec{t}, (\varphi,\psi)}(i)$ is given in
Eq.~(\mref{eq:coef-re-def1}) and
$\vec{a}\ssha_{(\varphi,\psi)}\vec{b}$ is given in
Eq.~(\mref{eq:mulind2}).
\end{theorem}

\begin{coro}
Let $\vec{r}\in \ZZ_{\geq 1}^k$ and $\vec{s}\in \ZZ_{\geq 1}^\ell$
with $r_1,s_1\geq 2$. Then
$$
\zeta(\vec{r})\,\zeta(\vec{s}) =   \sum_{\vec{t}\in \ZZ_{\geq
1}^{k+\ell}, |\vec{t}|=
  |\vec{r}|+|\vec{s}|}
  \Big(\sum _{(\varphi,\psi)\in \indI_{k,\ell}}
  \prod_{i=1}^{k+\ell}c_{\vec{r},\vec{s}}^{\vec{t},
(\varphi,\psi)}(i)\Big)
    \zeta(\vec{t})
$$
where $c_{\vec{r},\vec{s}}^{\vec{t}, (\varphi,\psi)}(i)$ is given in
Eq.~(\mref{eq:coef-re-def1}). \mlabel{co:mainmzv}
\end{coro}

\section{The algebraic framework of Connes and Kreimer on renormalization}
\mlabel{sec:ck}

The Algebraic Birkhoff Decomposition
%(henceforth abbreviated as \ABD)
of Connes and Kreimer is a fundamental result in their ground breaking work~\mcite{CK1} on Hopf algebra approach to renormalization of perturbative quantum field theory (pQFT).
This decomposition also links the physics theory of renormalization to Rota-Baxter algebra that has evolved in parallel to the development of QFT renormalization for several decades.

The introduction of Rota-Baxter algebra by G. Baxter~\mcite{Ba} in
1960 was motivated by Spitzer's identity~\mcite{Sp} that appeared in
1956 and was regarded as a
remarkable formula in the fluctuation theory of probability. Soon
Atkinson~\mcite{At} proved a simple yet useful factorization theorem
in Rota-Baxter algebras. The identity of Spitzer took its algebraic
form through the work of Cartier, Rota and Smith~\mcite{Ca,R-S}
(1972).

It was during the same period when the renormalization theory of pQFT was developed, through the the work of Bogoliubov and Parasiuk~\mcite{BP} (1957), Hepp~\mcite{He}(1966) and Zimmermann~\mcite{Zi} (1969), later known as the BPHZ prescription.

Recently QFT renormalization and Rota-Baxter algebra are tied together through the algebraic formulation of Connes and Kreimer for the former and a generalization of classical results on Rota-Baxter algebras in the latter~\mcite{EGK2,EGK3}. More precisely, generalizations of Spitzer's identity and Atkinson factorization give the twisted antipode formula and the algebraic Birkhoff decomposition in the work of Connes and Kreimer.

We recall the algebraic Birkhoff decomposition in Section~\mref{ss:ab}, prove the Atkinson factorization in Section~\mref{ss:af} and derive the algebraic Birkhoff decomposition from the Atkinson factorization in Section~\mref{ss:abaf}.

\subsection{Algebraic Birkhoff decomposition}
\mlabel{ss:ab}

For a $\bfk$-algebra $A$ and a $\bfk$-coalgebra $C$, we define the
convolution of two linear maps $f,g$ in $\Hom(C,A)$ to be
the map $f\star g\in\Hom(C,A)$ given by the composition
$$
C \xrightarrow{\Delta} C \ot C \xrightarrow{f\ot g} A \ot A
\xrightarrow{\mult} A.
$$
%In other words, $$(f*g)(a) = \sum_{(a)} f(a_{(1)}) \, g(a_{(2)}).$$
\begin{theorem} {\bf (Algebraic Birkhoff Decomposition)}
Let $H$ be a connected filtered Hopf algebra over $\CC$. Let $(A,\Pi)$ be a commutative Rota-Baxter algebra of weight $-1$ with $\Pi^2=\Pi$.
\begin{enumerate}
\item For
$\phi\in \mchar(H,A)$, there are unique linear maps $\phi_-:H\to \bfk+\Pi(A)$ and $\phi_+:H\to \bfk+(\id-\Pi)(A)$ such
that
\begin{equation}\label{eq:ABF}
 \phi=\phi_-^{\cprod (-1)}\cprod \phi_+.
\end{equation}
\mlabel{it:decom}
\item The elements $\phi_-$ and $\phi_+$
take the following forms on $\ker \vep$.
\begin{align}
\phi_-(x)&=
-\Pi(\phi(x)+\sum_{(x)}\phi_-(x')\phi(x'')),
\mlabel{eq:phi-}\\
\phi_+(x)&=
\tilde{\Pi}(\phi(x)+\sum_{(x)}\phi_-(x')\phi(x'')),
\mlabel{eq:phi+}
%\\ &= \tilde{Q}(\phi(x)+\sum_{(x)}\phi(x')\phi_+(x''))
\end{align}
\mlabel{it:rec}
where we have used the notation $\Delta(x)=1\ot x + x\ot 1 + \sum_{(x)} x'\ot x''$ with $x',x''\in \ker\vep$.
\item The linear maps $\phi_-$ and $\phi_+$
are also algebra homomorphisms. \mlabel{it:decgp}
\end{enumerate}
\mlabel{thm:algBi}
\end{theorem}
We call $\phi_+$ the {\bf renormalization} of $\phi$ and call $\phi_-$ the {\bf counter-term}.
Here is roughly how the renormalization method can be applied through the Algebraic Birkhoff Decomposition. See the tutorial article~\mcite{Gugn} for further details, examples and references.

Theorem~\mref{thm:algBi} can applied to renormalization as follows.
Suppose there is a set of divergent formal expressions, such as MZVs with not
necessarily positive arguments, that carries a certain algebraic combinatorial structure and from which we would like to extract finite values. On one hand, we first apply a suitable regularization (deformation) to each of these formal expressions so that the formal expression can be viewed as a singular value of the deformation function. Expanding around the singular point gives a Laurent series in $\bfk[\vep^{-1},\vep]]$. On the other hand, the algebraic combinatorial structure of the formal expressions, inherited by the deformation functions, can be abstracted to a free object in a suitable category. This free object parameterizes the deformation functions and often gives a Hopf algebra $H$. Thus the parametrization gives a morphism $\phi: H\to \bfk[\vep^{-1},\vep]]$ in the suitable category. Upon applying the Algebraic Birkhoff Decomposition, we obtain $\phi_+:H\to \bfk[[\vep]]$ which, composed with $\vep\mapsto 0$, gives us well-defined values in $\bfk$.

\subsection{Atkinson factorization}
\mlabel{ss:af}

The following is the classical result of Atkinson.
\begin{theorem}
{\rm (Atkinson Factorization)} Let $(R,P)$ be a Rota--Baxter algebra of weight $\lambda\neq 0$. Let $a\in R$. Assume that $b_\ell$ and $b_r$ are solutions of the fixed point equations
\begin{equation}
 b_\ell=1+P(b_\ell a), \qquad
b_r=1+(\id_R-P)(ab_r).
\mlabel{eq:recurlr}
\end{equation}
Then
$$
 b_\ell (1+\lambda a) b_r = 1.
$$
Thus
\begin{equation}
 1+\lambda a = b_\ell^{-1} b_r^{-1}
\mlabel{eq:atk3}
\end{equation}
if $b_\ell$ and $b_r$ are invertible.
\mlabel{thm:atk2} \mlabel{thm:Atkinson}
\end{theorem}

We note that the factorization (\mref{eq:atk3}) depends on the existence of invertible solutions of Eq.~(\mref{eq:recurlr}) that we will address next.

\begin{defn}
A {\bf filtered $\bfk$-algebra} is a $\bfk$-algebra $R$ together
with a decreasing filtration $R_n, \: n\geq 0, $ of nonunitary
subalgebras such that
$$
\bigcup_{n\geq 0} R_n = R, \quad R_n R_m \subseteq R_{n+m}.
$$
\end{defn}
It immediately follows that $R_0=R$ and each $R_n$ is an ideal of
$R$.
A filtered algebra is called {\bf complete} if $R$ is a complete
metric space with respect to the metric defined by the subsets $\{R_n\}$. Equivalently, a filtered $\bfk$-algebra $R$ with
$\{R_n\}$ is complete if $\cap_n R_n = 0$ and if the resulting
embedding
$$
R \to \bar{R}:= \varprojlim R/R_n
$$
of $R$ into the inverse limit is an isomorphism.

A Rota-Baxter algebra $(R,P)$ is called {\bf complete} if there are submodules $R_n\subseteq R, n\geq 0,$ such that $(R,R_n)$ is a complete algebra and $P(R_n)\subseteq R_n$.
\begin{theorem}
{\rm (Existence and uniqueness of the Atkinson factorization)}
Let $(R,P,R_n)$ be a complete Rota-Baxter algebra. Let $a$ be in $R_1$.
\begin{enumerate}
\item
The equations in (\mref{eq:recurlr}) have unique solutions $b_\ell$ and $b_r$. Further $b_\ell$ and $b_r$ are invertible. Hence Atkinson Factorization (\mref{eq:atk3}) exists.
\mlabel{it:atke}
\item
If $\lambda$ has no non-zero divisors in $R_1$ and $P^2=-\lambda P$ (in particular if $P^2=-\lambda P$ on $R$), then
there are unique $c_\ell\in 1+P(R)$ and $c_r\in 1+(\id_R-P)(R)$ such that
$$ 1+\lambda a= c_\ell c_r.$$
\mlabel{it:atku}
\end{enumerate}
\mlabel{thm:atke}
\end{theorem}

\subsection{From Atkinson factorization to algebraic Birkhoff decomposition}
\mlabel{ss:abaf}
We now derive the Algebraic Birkhoff Decomposition of Connes and Kreimer in Theorem~\mref{thm:algBi} from Atkinson Factorization in Theorem~\mref{thm:atke}.
Adapting the notations in Theorem~\mref{thm:algBi}, let $H$ be a connected filtered Hopf algebra and let $(A,Q)$ be a commutative Rota-Baxter algebra of weight $\lambda=-1$ with $Q^2=Q$, such as the pair $(A,Q)$ in Theorem~\mref{thm:algBi} (see also Example~\mref{ex:lau}). The increasing filtration on $H$ induces a decreasing filtration
$R_n=\{ f\in \Hom(H,A)\ |\ f(H^{n-1})=0 \}, n\geq 0$
on $R:=\Hom(H,A)$, making it a complete algebra. Further define
$$ P: R \to R, \quad P(f)(x)=Q(f(x)), f\in \Hom(H,A), x\in H.$$
Then it is easily checked that $P$ is a Rota-Baxter operator of weight $-1$ and $P^2=P$. Thus $(R,R_n,P)$ is a complete Rota-Baxter algebra.

Now let $\phi:H\to A$ be a character (that is, an algebra homomorphism). Consider $e-\phi:H\to A$. Then
$$ (e-\phi) (1_H) =e(1_H)-\phi(1_H)=1_H-1_H=0.$$
Thus $e-\phi$ is in $R_1$. Take $e-\phi$ to be our $a$ in Theorem~\mref{thm:atke}, we see that there are unique
$c_\ell\in P(R_1)$ and $c_r\in P(R_1)$ such that
$$ \phi = c_\ell c_r.$$
Further, by Theorem~\mref{thm:atk2}, for $b_\ell=c_\ell^{-1}$,
$b_\ell=e+P(b_\ell \cprod (e-\phi))$. Thus for $x\in \ker \vep =\ker e$, we have
\begin{eqnarray*}
b_\ell (x)&=& P(b_\ell \cprod (e-\phi))(x) \\
&=& \sum_{(x)} Q(b_\ell(a_{(1)})(e-\phi)(a_{(2)}))\\
&=& Q\big( b_\ell(1_H)(e-\phi)(x)+
\sum_{(a)} b_\ell(x')(e-\phi)(x'') + b_\ell(x)(e-\phi)(1_H)\\
&=& -Q\big( \phi(x) + \sum_{(x)} b_\ell (x')\phi(x'')\big).
\end{eqnarray*}
In the last equation we have used
$e(a)=0, e(a'')=0$ by definition.
Since $b_\ell(1_H)=1_H$, we see that $b_\ell=\phi_-$ in Eq.~(\mref{eq:phi-}).

Further, we have
$$c_r=c_\ell^{-1} \phi=b_\ell \phi
= -b_\ell(e-\phi) + b_\ell
= -b_\ell(e-\phi) + e+P(b_\ell (e-\phi))
= e-(\id-P)(b_\ell (e-\phi)).$$
With the same computation as for $b_\ell$ above,
we see that $c_r = \phi_+$ in Eq.~(\mref{eq:phi+}).

\section {Heat-kernel type regularization approach to the renormalization of MZVs}
\mlabel{sec:gz}

To extend the double shuffle relations to MZVs with non-positive arguments, we
have to make sense of the divergent sums defining these MZVs. For this
purpose, we adapt the renormalization method from quantum field theory in the
algebraic framework of Connes-Kreimer recalled in the last section.
We will give three approaches including the approach in this section using a heat-kernel type regularization, named after a similar process in physics. Since examples and motivations of this approach can already be found elsewhere~\mcite{Gusu,GZ,GZ2}, we will be quite sketchy in this section. More details will be given to the two other approaches in Section~\mref{sec:mp}.

\subsection{Renormalization of MZVs}
Consider the abelian semigroup
\begin{equation}
\frakM= \{{{\wvec{s}{r}}}\ \big|\ (s,r)\in \ZZ \times \RR_{>0}\}
\end{equation}
with the multiplication
$$ {\wvec{s}{r}}\cdot {\wvec{s'}{r'}}={\wvec{s+s'}{r+r'}}.$$
With the notation in Section~\mref{ss:msh}, we define the Hopf algebra
$$\calh_{\frakM}:=\sh_{\CC,1}(\CC \frakM)$$
with the quasi-shuffle product $*$ and the deconcatenation coproduct $\Delta$ in Section~\mref{ss:msh}.
For $w_i=\wvec{s_i}{r_i}\in \frakM,\ i=1,\cdots,k$, we use the notations
$$ \vec{w}=(w_1, \dots,w_k)
=\wvec{s_1, \dots,s_n}{r_1, \dots,r_k}=\wvec{\vec s}{\vec r},\ {\rm
where\ } \vec s=(s_1, \dots,s_k), \vec r=(r_1, \dots,r_k).$$

For $\vec{w}=\wvec{\vec s}{\vec r}\in \frakM^k$ and $\vep\in \CC$
with ${\rm Re}(\vep)<0$, define the {\bf directional regularized
MZV}:
\begin{equation}
Z(\wvec{\vs}{\vr};\vep)=\sum_{n_1>\cdots>n_k>0}
\frac{e^{n_1\,r_1\vep} \cdots
    e^{n_k\,r_k\vep}}{n_1^{s_1}\cdots n_k^{s_k}}
\label{eq:reggmzv}
\end{equation}
It converges for any $\wvec{\vs}{\vr}$ and is regarded as the
regularization of the {\bf formal MZV}
\begin{equation}
\zeta (\vs)= \sum_{n_1>\cdots>n_k>0} \frac{1}{n_1^{s_1} \cdots
    n_k^{s_k}}
\label{eq:formgmzv}
\end{equation}
which converges only when $s_i>0$ and $s_1>1$. It is related to the
multiple polylogarithm
$${\rm Li}_{s_1, \dots,s_k}(z_1, \dots,z_k)=\sum_{n_1>\cdots > n_k>0}
    \frac{z_1^{n_1} \cdots z_k^{n_k}}{n_1^{s_1}\cdots n_k^{s_k}}$$
by a change of variables $z_i=e^{r_i\vep}, 1\leq i\leq k$.

This regularization defines an algebra homomorphism~\mcite{GZ}:
\begin{equation}
\uni {Z}: \calh_\frakM \to \CC[T][[\vep,\vep^{-1}], \label{eq:zmap}
\end{equation}
In the same way, for
\begin{equation}
\frakM ^-= \big\{{{\wvec{s}{r}}}\ \big|\ (s,r)\in \ZZ _{\le 0} \times
\RR_{>0}\big\},
\end{equation}
$\uni{Z}$ restricts to an algebra homomorphism
\begin{equation}
\uni{Z}: \calh_{\frakM^-} \to R:=\CC[[\vep,\vep^{-1}].
\mlabel{eq:zmap-}
\end{equation}

Since both $(\CC[T][\vep^{-1},\vep]],\Pi)$ and $(\CC[\vep^{-1},\vep]],\Pi)$, with $\Pi$ defined in Example~\mref{ex:lau}, are commutative Rota-Baxter algebras with $\Pi^2=\Pi$,
 we have the decomposition
$$ \uni{Z}=\uni{Z}_-^{-1} \star \uni{Z}_+$$
by the algebraic Birkhoff decomposition in Theorem~\mref{thm:algBi} and obtain

\begin{theorem} (\mcite{GZ,GZ2})
The map $\uni{Z}_+: \calh_\frakM\to \CC[T][[\vep]]$ is an algebra
homomorphism which restricts to an algebra homomorphism $\uni{Z}_+:
\calh_{\frakM^-}\to \CC[[\vep]]$. \label{thm:renormz}
\end{theorem}

Because of Theorem~\mref{thm:renormz}, the following definition is
valid.
\begin{defn}
For $\vs=(s_1, \dots,s_k)\in \ZZ^k$ and $\vr=(r_1, \dots,r_k)\in
\RR_{>0}^k$, define the {\bf renormalized directional MZV} by
\begin{equation} \zeta\lp \wvec{\vs}{\vr}\rp = \lim_{\vep\to 0} \uni{Z}_+\lp \wvec{\vs}{\vr};\vep\rp .
\label{eq:dmzv}
\end{equation}
Here $\vr$ is called the {\bf direction vector}. \label{de:dmzv}
\end{defn}

As a consequence of Theorem~\mref{thm:renormz}, we have
\begin{coro}
The renormalized directional MZVs satisfy the quasi-shuffle
relation
\begin{equation}
\zeta\lp \wvec{\vs}{\vr}\rp \zeta\lp \wvec{\vs\,'}{\vr\,'}\rp =
\zeta\lp \wvec{\vs}{\vr} \msh \wvec{\vs\,'}{\vr\,'} \rp.
\label{eq:dqshuf}
\end{equation}
Here the right hand side is defined in the same way as in
Eq.~(\mref{eq:quasi}). \label{co:qshuf}
\end{coro}

\begin{defn}
For $\vec{s}\in \ZZ_{> 0}^k\cup \ZZ_{\leq 0}^k$,  define
\begin{equation} \zeta\lp \vec{s}\rp
= \lim_{\delta \to 0^+} \zeta\lp \wvec{\vec{s}}{|\vec{s}|+\delta}\rp
, \mlabel{eq:gmzv}
\end{equation}
where, for $\vec s=(s_1,\cdots,s_k)$ and $\delta\in \RR_{>0}$, we
denote $|\vec{s}|=(|s_1|,\cdots,|s_k|)$ and $|\vec
s|+\delta=(|s_1|+\delta,\cdots,|s_k|+\delta).$ These
$\zeta(\vec{s})$ are called the {\bf renormalized MZVs} of the
multiple zeta function $\zeta(u_1,\cdots,u_k)$ at $\vec{s}$.
\mlabel{de:rmzv}
\end{defn}

\begin{theorem} {\bf \mcite{GZ}}
\begin{enumerate}
\item
The limit in Eq.~(\mref{eq:gmzv}) exists for any $\vec{s}=(s_1,\cdots,s_k)\in \ZZ_{> 0}^k\cup \ZZ_{\leq 0}^k$.
\item
When $s_i$ are all positive with $s_1>1$, we have $\zeta\lp \wvec{\vec{s}}{\vec{r}}\rp =\zeta(\vec{s})$ independent of $\vec{r}\in \ZZ_{>0}^k$. In particular, we have $\gzeta(\vec{s})=\zeta(\vec{s})$.
\mlabel{it:g1}
\item
When $s_i$ are all positive, we have $\gzeta(\vec{s})=\zeta\lp \wvec{\vec{s}}{\vec{s}}\rp $. Further, $\gzeta(\vec{s})$ agrees with the regularized MZV $Z^*_{\vec{s}}(T)$ defined by Ihara-Kaneko-Zagier~\mcite{IKZ}.
\mlabel{it:ge1}
\item
When $s_i$ are all negative, we have
$\gzeta(\vec{s})=\zeta\lp \wvec{\vec{s}}{-\vec{s}}\rp =\disp{\lim_{\vec{r}\to -\vec{s}} \zeta\lp \wvec{\vec{s}}{\vec{r}}\rp}.$ Further, these values are rational numbers.
\mlabel{it:neg}
\item
The value $\gzeta(\vec{s})$ agrees with $\zeta(\vec{s})$ whenever the latter is defined by analytic continuation.
\mlabel{it:nonp}\\
\noindent
{\ ${}$ } \hspace{-1.70cm} Furthermore,
\item
the set $\{\gzeta(\vec{s})\big| \vec{s}\in \ZZ^k_{>0}\}$ satisfies the quasi-shuffle relation;
\mlabel{it:pqs}
\item
the set $\{\gzeta(\vec{s})\big| \vec{s}\in \ZZ^k_{\leq 0}\}$ satisfies the quasi-shuffle relation.
\mlabel{it:nqs}
\end{enumerate}
\mlabel{thm:main}
\end{theorem}

The following table lists $\gzeta(-a,-b)$ for $0\leq a,
b\leq 6$.

{\small
\begin{equation}
\label{eq:gzt}
%{\ } \vspace{10.5cm} {\ }\vspace{3cm}
%\hspace{1.5cm}
\begin {array} {c|c|c|c|c|c|c|c|}
\gzeta(-a,-b)&a=0&a=1 &a=2&a=3&a=4&a=5&a=6\\ \hline
&&&&&&&\\
b=0 & \frac 38&\frac 1{12}&\frac 1{120}&-\frac 1{120}&-\frac
1{252}&\frac 1{252}&\frac 1{240}
\\
&&&&&&&\\\hline
&&&&&&&\\
b=1 & \frac 1{24}&\frac 1{288}& -\frac 1{240}& \frac {83}{64512}&
\frac 1{504}& -\frac {3925}{2239488}& -\frac 1{480}
\\
&&&&&&&\\\hline
&&&&&&&\\
b=2 & -\frac 1{120}&-\frac 1{240}& 0& \frac 1{504}& -\frac
{319}{437400}&
-\frac 1{480}& \frac {2494519}{1362493440}\\
&&&&&&&\\\hline
&&&&&&&\\
b=3 & -\frac 1{240}&-\frac {71}{35840}& \frac 1{504}& \frac
1{28800}& -\frac 1{480}& \frac {114139507}{139519328256}& \frac
1{264}
\\
&&&&&&&\\\hline
&&&&&&&\\
b=4 & \frac 1{252}&\frac 1{504}& \frac {319}{437400}& -\frac 1{480}&
0& \frac 1{264}& -\frac {41796929201}{26873437500000}
\\
&&&&&&&\\\hline
&&&&&&&\\
b=5 &\frac 1{504}& \frac {32659}{15676416}& -\frac 1{480}& -\frac
{21991341}{25836912640}& \frac 1{264}& \frac 1{127008}& -\frac
{691}{65520}
\\
&&&&&&&\\\hline
&&&&&&&\\
b=6 &-\frac 1{240} &-\frac 1{480}& -\frac {2494519}{1362493440}&
\frac 1{264}& \frac {41796929201}{26873437500000}&
-\frac {691}{65520}& 0\\
&&&&&&&\\\hline
\end {array}
\end{equation}
}

\subsection{The differential structure}
The shuffle relation for convergent MZVs from their integral representations does not directly generalize to renormalized MZVs due to the lack of a suitable integral representation. However a differential variation of the shuffle relation might exist for renormalized MZVs. One evidence is the following differential version of the algebraic Birkhoff decomposition~\mcite{GZ2} for renormalized MZVs and further progress will be discussed in a paper under preparation. We first recall some concepts.

\begin{enumerate}
\item
A {\bf differential algebra} is a pair $(A,d)$ where $A$ is an algebra and $d$ is a {\bf differential operator}, that is, such that  $d(xy)=d(x)y+xd(y)$ for all $x,y\in A$. A differential algebra homomorphism $f: (A_1,d_1)\to (A_2,d_2)$ between two differential algebras $(A_1,d_1)$ and $(A_2,d_2)$ is an algebra homomorphism $f:A_1\to A_2$ such that $f\circ d_1 = d_2 \circ f$.
\item
A {\bf differential Hopf algebra} is a pair $(H,d)$ where $H$ is a
Hopf algebra and $d:H\to H$ is a differential operator such that
\begin{equation}
\Delta(d(x))=
\sum_{(x)} \big( d(x_{(1)})\bigotimes\, x_{(2)} +
    x_{(1)} \bigotimes\, d(x_{(2)})\big).
\mlabel{eq:diffH0}
\end{equation}
\item
A {\bf differential Rota-Baxter algebra} is a triple $(A,\Pi,d)$ where
$(A,\Pi)$ is a Rota-Baxter algebra and $d:R\to R$ is a differential
operator such that $P\circ d=d\circ P$.
\end{enumerate}

\begin{theorem}
{\bf (Differential Algebraic Birkhoff Decomposition)~\mcite{GZ2}} Under the same assumption as in Theorem~\mref{thm:algBi}, if in addition $(H,d)$ is a differential Hopf algebra, $(A,\Pi,\partial)$ is
a commutative differential Rota-Baxter algebra,
and $\phi:H\to A$ is a differential algebra homomorphism,
then the maps $\phi_-$ and $\phi_+$ in Theorem~\mref{thm:algBi} are also differential algebra homomorphisms.
\mlabel{thm:dabd}
\end{theorem}

\begin{theorem} {\bf (\mcite{GZ2})}
\begin{enumerate}
\item
For $\wvec{s}{r}\in \frakM$, define
$d(\wvec{s}{r}=r\wvec{s-1}{r}$. Extend $d$ to $\calh_\frakM
=\oplus_{k\geq 0} (\bfk \frakM)^{\ot k}$ by defining, for $\fraka:=a_1\ot \cdots\ot a_k\in (\bfk \frakM)^{\ot k}$,
\begin{equation}
 d(\fraka)=\sum_{i=1}^k a_{i,1}\ot\cdots\ot a_{i,k},
 \quad a_{i,j}=\left \{ \begin{array}{ll} a_j, & j\neq i, \\
    d(a_j),& j=i. \end{array} \right .
 \mlabel{eq:diffH}
 \end{equation}
Then $(\calh_A,d)$ is a differential Hopf algebra.
\item
The triple $(\CC[\vep^{-1},\vep]],\Pi,\frac{d}{d\vep})$ is a commutative differential Rota-Baxter algebra.
\item
The map
$\uni {Z}: \calh_\frakM \to \CC[[\vep,\vep^{-1}] $ defined in Eq.~(\mref{eq:zmap-}) is a differential algebra homomorphism.
\item
The algebra homomorphism $\uni{Z}_+:
\calh_{\frakM^-}\to \CC[[\vep]]$ in Theorem~\mref{thm:renormz} is a differential algebra homomorphism.
\end{enumerate}
\mlabel{thm:dmzv}
\end{theorem}

%%%%%%%%%%%%%%%%%%%%%%%%%%%%%%%%%%%%%%%%%%%%%%%%%%%%%%%
%%%%%%%%%%%%%%%%%%%%%%%%%%%%%%%%%%%%%%%%%%%%%%%%%%%%%%

\section{Renormalization of multiple zeta values seen as nested  sums of symbols}
\mlabel{sec:mp}

We present two more approaches to renormalize multiple zeta functions at non-positive integers,  both of which  lead to MZVs which obey stuffle
relations. Like the renormalization method described in the previous
  section, they both give rise to rational multiple zeta values at
non-positive integers and we check that  the two methods yield the same double
multiple zeta values at non-positive integer arguments.  This presentation is
based on joint work of one of the authors with D. Manchon
\mcite{MP2} in which  multiple zeta functions are viewed as particular
  instances of nested sums of symbols and where the algebraic Birkhoff
decomposition approach is used to renormalize multiple zeta
functions at poles. Here, we furthermore present an alternative renormalization method
based on generalized  evaluators used in physics \cite{S}.

\subsection{A class of symbols}
For a complex number  $b$,  a smooth function $f:\RR -\{0\}\to \CC$ is called {\bf positively homogeneous of degree $b$} if
$f(t\,  \xi)=t^b f(\xi)$ for all $t> 0$ and $\xi\in \RR$.

The symbols which were originally defined on $\R^n$ are now defined
on $\R$ which is sufficient for our needs in this paper. We call a smooth function $\sigma:\RR\to \CC$ a {\bf symbol} if there is a
  real number $a$ such that
for any non-negative integer $\gamma$, there is a positive constant $C_{\gamma}$ with
$$\left\vert \partial^\gamma \sigma(\xi)\right\vert \le C_{\gamma}(1+\vert
\xi\vert)^{a-\gamma}, \quad \forall \xi\in
\RR.$$
For a complex number $a$ and a  non-negative integer  $j$, let
$\sigma_{a-j}: \RR-\{0\}\to \C$  be  a smooth and positively homogeneous function of degree
$a-j$. We write
$\sigma\sim \sum_{j=0}^\infty  \sigma_{a-j}$ if, for any non-negative integer $N$ and non-negative integer $\gamma$, there is a positive constant $C_{\gamma,N}$ such that
$$\left\vert \partial^\gamma\left(\sigma(\xi)-\sum_{j=0}^{N}\,
    \sigma_{a-j}(\xi)\right)\right\vert \le C_{\gamma,N}(1+\vert
\xi\vert)^{\smop{Re }(a)-N-1-\gamma}, \quad \forall \xi \in \RR-\{0\},$$
 where ${\rm Re }(a)$ stands for the real part of $a$.

For any complex number  $a$ and any non-negative integer $k$, a  symbol  $\sigma:\RR\to \CC$ is called a {\bf  log-polyhomogeneous  of log-type $k$ and order $a$} if
\begin{equation}\mlabel{eq:asymptf}
\sigma(\xi)=\sum_{l=0}^k \sigma_l(\xi)\log ^l|\xi|, \quad
\sigma_l(\xi)\sim\sum_{j=0}^\infty \sigma_{a-j, l}(\xi)
\end{equation}
with  $\sigma_{a-j, l}(\xi)$  positively homogeneous of degree $a-j$.

Let
${\mathcal S}^{a, k}$ denote the linear space over $\CC$
of log-polyhomogeneous symbols on $\R$
of log-type $k$  and order $a$.
Then we have $\cals^{a,k}\subseteq \cals^{a,k+1}$.
Let ${\mathcal S}^{*, k}$ denote the linear span over $\C$ of all ${\mathcal S}^{a,k}$ for
${a\in \C}$. Then ${\mathcal S}^{*,0}$ corresponds to the algebra of {\rm classical symbols} on $\R$.
  We also define
$${\mathcal S}^{*, *}:=
\bigcup_{k=0}^\infty {\mathcal S}^{*, k}$$
which is an  algebra for
the ordinary product of functions  filtered by the log-type \mcite{L}
  since the product  of two symbols of log-types $k$ and
$k^\prime$ respectively is of log-type $k+k^\prime$.
 The union
$\bigcup_{a\in \Z}\bigcup_{k=0}^\infty {\mathcal S}^{a, k}$ is a
  subalgebra of ${\mathcal S}^{*,*}$, and  $\bigcup_{a\in \Z} {\mathcal S}^{a, 0}$ is a
  subalgebra of ${\mathcal S}^{*,0}$.

Let  ${\mathcal
  P}^{\alpha,k}$ be the algebra of  {\bf positively supported\/} symbols,
i.e. symbols in ${\mathcal S}^{\alpha,k}$ with support  in
$(0,+\infty)$ so that they are non-zero only at positive arguments.
We keep {\sl mutatis mutandis\/} the above notations; in
particular ${\mathcal P}^{*,0}$ is a subalgebra of the filtered algebra  ${\mathcal
  P}^{*,*}$.

For $\sigma \in   {\mathcal
  P}^{\alpha,k}$  we call $\mopl{fp}_{\xi\to \infty} \sigma(\xi):=\sigma_{0,0}(\xi)$
the {\bf finite part}  at zero (so named since it it reminiscent of
  Hadamard's finite parts) of such a symbol  $\sigma$ which corresponds to the
  constant term in the expansion.

The following rather elementary
  statement is our  main motivation  here for  introducing log-polyhomogeneous
  symbols.
\begin{prop}\mlabel{prop:contRotaBaxter} \mcite{MP1} The operator $I$
defined in (\mref{eq:int})  on the algebra $C[0,\infty)$ by
$$f\mapsto \left(\xi\mapsto  I(f)(\xi)=\int_0^\xi f(t)dt\right)
$$
maps ${\mathcal P}^{*, k-1}$ to ${\mathcal P}^{*, k}$ for any positive integer $k$.
\end{prop}
By Proposition \mref{prop:contRotaBaxter}, for any $\sigma$ in ${\mathcal
  P}^{*,k}$, the primitive $I(\sigma)(\xi)$  has an asymptotic behavior as
$\xi\to \infty$  of the
type (\mref{eq:asymptf}) with $k$ replaced by $k+1$. The constant term defines
the cut-off regularized integral (see e.g. \mcite{L}):
$$
\cutoffint_0^\infty \sigma(t)\, dt:= \mopl{fp }_{\xi \to \infty} \int_0^\xi
\sigma(t)\, dt.
$$

\subsection{Nested sums of symbols and their pole structures}
\mlabel{ss:nsp}
\subsubsection{Nested sums}
Recall that the operator $I $ on ${\mathcal P}^{*,*}$ defined by Eq.~(\mref{eq:int})
satisfies the weight zero Rota-Baxter relation (\mref{eq:Ba}).
On the other hand the operator $P$ defined by Eq.~(\mref{operatorP})
satisfies the Rota-Baxter relation with weight $\lambda =-1$ and the operator
$Q=P-Id$ in Eq.~(\mref{operatorQ})
satisfies the Rota-Baxter relation with weight $\lambda =1$.

The Rota-Baxter operators  $P$ and $I$ relate by  means of the Euler-MacLaurin formula which compares discrete sums with
integrals.
For
$\sigma\in{\mathcal P}^{*,*}$  the
Euler-MacLaurin formula
(see e.g. \mcite{Ha})  reads:
\begin{eqnarray}\label{eq:EML}
P(\sigma)(N)-I(\sigma)(N)&=&\frac 12 \sigma(N)+\sum_{k=2}^{2K}\frac{B_{k}}{k!}\sigma^{(k-1)}(N)\nonumber\\
&+&
\frac{1}{(2K+1)!}\int_0^N\overline{B_{2K+1}}(x)\sigma^{(2K+1)}(x)\,dx.
\end{eqnarray}
with $\overline{B_k}(x)= B_k\left(x-[x] \right)$. Here $
B_k(x)= \sum_{i=0}^k {k\choose i} \, B_{k-i} \, x^k$ are the Bernoulli
polynomials of degree $k$,
the $B_i$ being the Bernoulli numbers, defined by the generating series:
$$\frac {t}{e^t-1}=\sum_i\frac{B_i}{i!}t^i.$$
Since
$B_k(1)=B_k$ for any  $k\geq 2$, setting
$x=1$ we have \begin{equation}\mlabel{eq:sumBi}
B_k= \sum_{i=0}^k {k\choose i} \, B_{k-i} =\sum_{i=0}^k {k\choose i} \,
B_{i},\quad \forall k\geq 2.
\end{equation}

The Euler-MacLaurin formula therefore provides an interpolation of $P(\sigma)$ by a symbol.
\begin{prop} \mlabel{prop:interpolation}\mcite{MP2}
For any $\sigma\in {\mathcal P}^{a, k}$, the discrete sum $P(\sigma)$ can be interpolated by
a symbol  $\overline P(\sigma)$ in ${\mathcal P}^{a+1, k+1}+{\mathcal
  P}^{0, k+1} $   (i.e. $\overline P(\sigma)(n)= P(\sigma)(n)=\sum_{k=0}^n\sigma(k),\quad \forall n\in \N$) such that  $$\overline
  P(\sigma)- I(\sigma) \in
  {\mathcal P}^{a, k}.$$
The operator $\overline Q:= \overline P-Id: {\mathcal P}^{a, k}\to {\mathcal P}^{a+1, k+1}+{\mathcal
  P}^{0, k+1}$  interpolates $Q$.
\end{prop}
By Proposition \mref{prop:interpolation}, given  a symbol $\sigma$ in $ {\mathcal P}^{a, k}$, the
interpolating  symbol  $\overline P(\sigma)$ lies in  ${\mathcal P}^{a+1, k+1}+{\mathcal
  P}^{0, k+1} $.   It follows that the discrete sum  $P(\sigma)(N)=\overline {
  P}(\sigma)(N)$ has  an  asymptotic behavior for large $N$
given by finite linear combinations of expressions of the type
(\mref{eq:asymptf}) with $k$ replaced by $k+1$ and $a$ by $a+1$ or $0$.
Picking the finite part, for any $\sigma\in {\mathcal P}^{*, *}$  we define the following cut-off sum:
 \begin{equation}\mlabel{eq:cutoffsum}\cutoffsum_0^\infty  \sigma:= \mopl{fp }_{N\to \infty} P(\sigma)(N)=\mopl{fp }_{N\to \infty}
 \sum_{k=0}^N\sigma(k),\end{equation} which extends the ordinary discrete sum $\sum_0^\infty$ on
 $L^1$-symbols.
If $\sigma$ has non-integer order, we have $\cutoffsum_0^\infty  \sigma= \mopl{fp }_{N\to \infty}
 \sum_{k=0}^{N+K}\sigma(k)$ for any  integer $K$, so that in particular
$\cutoffsum_0^\infty \sigma=  \mopl{fp }_{N\to \infty} Q(\sigma)(N)$
since the operators $P$ and $Q$ only differ by an integer in the
  upper bound of the sum.

With the help of  the
 interpolation map  described in
 Proposition \mref{prop:interpolation}, we can assign to a tensor product
 $\sigma:=\sigma_1\otimes\cdots\otimes \sigma_k$  of (positively supported)
 classical symbols, two log-polyhomogeneous symbols defined inductively in the
 degree $k$ of the tensor product, which interpolate the nested
 iterated sum $$ \sum_{0\leq n_k\leq n_{k-1}\leq \cdots \leq n_{2}\leq
   n_1} \sigma_1(n_1)\cdots \sigma_k(n_k)=\sigma_1\, P\Big(\cdots \sigma_{k-2}\, P\big(\sigma_{k-1}\,
 P(\sigma_k)\big)...\Big),$$
$$ \sum_{0\leq n_k< n_{k-1}< \cdots < n_{2}< n_1} \sigma_1(n_1)\cdots \sigma_k(n_k)=\sigma_1\, Q\Big(\cdots \sigma_{k-2}\, Q\big(\sigma_{k-1}\,
 P(\sigma_k)\big)...\Big).$$
In the following we will only consider the second class of symbols, including
their regularization, renormalization and application to multiple zeta
values. A parallel approach applies to the first class of symbols with application to non-strict multiple zeta values in Eq.~(\mref{eq:mszv})~\mcite{OZ,Zud}.

\begin{thm} \mlabel{thm:sigmatilde}\mcite{MP2}
Given $\sigma_i\in{\mathcal P}^{\alpha_i,0}$, $i=1,\ldots ,k$, setting
$
\sigma:=\sigma_1\otimes\cdots\otimes \sigma_k$, the function $\widetilde \sigma$ defined by:
\begin{equation}\mlabel{eq:sigmatilde}
{\widetilde \sigma}:=\sigma_1\, \overline Q\Big(\cdots \sigma_{k-2}\,\overline Q\big(\sigma_{k-1}\,
 \overline Q(\sigma_k)\big)...\Big)
\end{equation}
which interpolates  nested sums in the following way:
\begin{eqnarray*}
 \quad\widetilde \sigma (n_1)&=& \sum_{0\leq n_k< n_{k-1}< \cdots
  <n_{2}< n_1} \sigma_1(n_1)\cdots
\sigma_k(n_k), \quad \forall n_1\in \N,
\end{eqnarray*}
lies in ${\mathcal P}^{*, k-1}$ as  linear combinations of
(positively supported) symbols in  ${\mathcal P}^{\alpha_1+\cdots +\alpha_j+j-1, j-1}$,
$j\in\{1,\ldots ,k\}$.
\end{thm} On the grounds of this result, we  define the cut-off nested discrete sum
of a tensor product of (positively supported) classical symbols.
\begin{defn}
{\rm
For $\sigma_1, \ldots, \sigma_k \in {\mathcal P}^{*,0}$ and
  $\sigma:=\sigma_1\otimes\cdots\otimes \sigma_k$
we call
\begin{equation*}
\cutoffsum^{\smop{Chen}}_{<}\sigma
:=\cutoffsum_{0}^\infty \widetilde \sigma
= \cutoffsum_{0<  n_k<\cdots < n_1} \sigma_1(n_1) \cdots  \sigma_k(n_k)
\end{equation*}
the {\bf cut-off nested sum} of $f=\sigma_1\otimes \cdots \otimes \sigma_k$.
}
\end{defn}

\subsubsection{The pole structure of nested sums of symbols}
To build meromorphic extensions, we combine the cut-off sum
$\cutoffsum_0^\infty$ introduced in (\mref{eq:cutoffsum}) with holomorphic
deformations of the symbol in the integrand.

A family $\{a(z)\}_{z\in \Omega}$  in a topological vector space ${\mathcal A}$ which is
parameterized by a complex domain $\Omega$, is holomorphic at
$z_0\in \Omega$ if
the corresponding function $f:\Omega\to{\mathcal A}$ admits a  Taylor expansion in a neighborhood $N_{z_0}$ of $z_0$
$$
a(z) = \sum_{k=0}^{\infty}a^{(k)}(z_0)\,\frac{(z-z_0)^k}{k!}
$$
which is convergent, uniformly on compact subsets of $N_{z_0}$ (i.e. locally uniformly),
with respect to the  topology on ${\mathcal A}$. The vector spaces of functions we consider here are $C(\R,\C)$ and
$C^\infty(\R,\C)$ equipped with their usual topologies, namely uniform
convergence on compact subsets, and uniform convergence of all derivatives on
compact subsets respectively.
\begin{defn}\mlabel{defn:holfamilies}
{\rm
Let $k$ be a non-negative integer, and
  let $\Omega$ be a domain in $\C$. A {\bf simple holomorphic
family of log-polyhomogeneous symbols\/} $\sigma(z)\in {\mathcal
S}^{*, k}$ parameterized by $\Omega$  is  a holomorphic family
$\sigma(z)(\xi) := \sigma(z,\xi)$  of smooth functions on $\R$ such that:
\begin{enumerate}
\item  the  order
$\alpha:\Omega \to\C$ is holomorphic on $\Omega$,
\item $\sigma(z)(\xi)=\sum_{l=0}^k \sigma_l(\xi)\, \log^l \vert \xi\vert$
with $$
 \sigma_l(z)(\xi) \sim \sum_{j\geq 0}\,
 \sigma(z)_{\alpha(z)-j,l}(\xi).$$
Here  $\sigma(z)_{\alpha(z)-j,l}$ positively homogeneous of degree $\alpha(z)-j$,
\item for any positive integer $N$ there is some positive integer $K_N$ such
  that the remainder term
$$\sigma_{(N)}(z)(\xi):= \sigma(z)(\xi) -\sum_{l=0}^k \sum_{j= 0}^{K_N}
 \sigma(z)_{\alpha(z)-j, l}(\xi)\log^l\vert\xi\vert=o(\vert \xi\vert^{-N})$$
is holomorphic in $z\in \Omega$ as a function of $\xi$ and verifies for any
$\epsilon>0$ the following estimates:
$$
\partial_\xi^\beta\partial_z^k\sigma_{(N)}(z)(\xi)=o(|\xi|^{-N-|\beta|+\e})
$$
locally uniformly in $z\in \Omega$ for $k\in \N$ and $\beta\in\N^n$.
\end{enumerate}
A {\bf holomorphic
family of log-polyhomogeneous symbols\/} is a finite linear combination (over
$\C$) of simple holomorphic families.
}
\end{defn}
It follows from the Euler-MacLaurin formula (see e.g. \cite{Ha,MP2}) that for any holomorphic family $\sigma(z)$
of symbols in ${\mathcal P}^{*,*}$, we have
$$\cutoffsum_{n=0}^\infty \sigma(z)(n)= \cutoffint_{0}^\infty \sigma(z)(\xi)\, d\xi
+C(\sigma(z))$$
with $z\mapsto C(\sigma(z))$  a holomorphic function at zero.
Hence,
$z\mapsto \cutoffsum_{n=0}^\infty \sigma(z)(n)$ and $z\mapsto
\cutoffint_{0}^\infty \sigma(z)(\xi)\, d\xi$ have the same pole structure.   Results by Kontsevich and Vishik \mcite{KV} for classical symbols and their generalization
by Lesch \mcite{L} to log-polyhomogeneous symbols, and relative to  the
pole structure  of
cut-off integrals of holomorphic families of symbols,
  therefore carry out to
discrete  cut-off sums of holomorphic (positively supported)
log-polyhomogeneous symbols.
Let us briefly recall the notion of holomorphic regularization
inspired by \mcite{KV}.
\begin{defn}
{\rm
A  {\bf holomorphic
regularization procedure} on ${\mathcal S}^{*, *}$
is a  map
\begin{eqnarray*}
{\mathcal R}: {\mathcal S}^{*, *}&\to & {\rm Hol}_\Omega\,\left( {\mathcal S}^{*, *}\right)\\
f  &\mapsto & \{\sigma(z)=\sigma_f(z)\}_{z\in \Omega}
\end{eqnarray*}
where $\Omega$ is an open subset of $\C$ containing $0$, and ${\rm Hol}_\Omega\left( {\mathcal S}^{*, *}\right)$ is the algebra of
holomorphic families in ${\mathcal S}^{*, *}$ ,
such that for any $f\in {\mathcal S}^{*, *}$,
\begin{enumerate}
\item  $\sigma(0)=f$,
\item  the holomorphic family $\sigma(z)$ can be written as a linear combination of
  simple ones:
$$
\sigma(z)=\sum_{j=1}^k\sigma_j(z),
$$
the holomorphic  order  $\alpha_j(z)$ of which verifies $\mop{Re
  }(\alpha_j^\prime(z))< 0$ for any $z\in \Omega$ and any $j\in\{1,\ldots ,k\}$.\\
\end{enumerate}
A  holomorphic regularization ${\mathcal R}$ is {\bf simple} if, for any
  log-polyhomogeneous symbol $\sigma\in{\mathcal S}^{\alpha,k}$, the holomorphic family
  ${\mathcal R}(\sigma)$
  is simple. Since we only consider simple holomorphic
    regularizations, we drop the explicit mention of simplicity.
}
\end{defn}
A similar definition holds with suitable subalgebras of ${\mathcal S}^{*,*}$,
e.g. classical symbols ${\mathcal S}^{*,0}$ instead of
log-polyhomogeneous. Holomorphic regularization procedures naturally arise in physics:
\begin{ex}\mlabel{ex:regphys}
{\rm  Let $z\mapsto \tau(z)  \in {\mathcal S}^{*,0}$ be a holomorphic family of
  classical symbols such that $\tau(0)=1$ and $\tau(z)$ has holomorphic  order
 $\alpha(z)$ with $\mop{Re }(\alpha^\prime(z))< 0$. Then
$${\mathcal R}: \sigma \mapsto \sigma(z):= \sigma\, \tau(z)$$
yields a holomorphic
regularization on ${\mathcal S}^{*, *}$ as well as on ${\mathcal
S}^{*, 0}$. Choosing $\tau(z)(\xi):=
\chi(\xi)+\big(1-\chi(\xi)\big)\big(H(z)\, \vert\xi \vert^{-z}\big)$
  where $H$ is a scalar valued  holomorphic map such that  $H(0)=1$, and where
  $\chi$ is a smooth cut-off function which is identically one
    outside the unit interval and zero in a small neighborhood of zero, we get
$$
{\mathcal R}(\sigma)(z)(\xi)= \chi(\xi)\sigma(\xi)+\big(1-\chi(\xi)\big)\big(H(z)\, \sigma(\xi)\, \vert \xi \vert^{-z}\big).
$$
Dimensional regularization commonly used in physics is of this type, where $H$
is expressed in terms of Gamma functions which account for a
``complexified'' volume of the unit sphere. When $H\equiv 1$, such a
regularization   ${\mathcal
  R}$ is called   Riesz regularization.
}
\end{ex}
\begin{prop}\mlabel{prop:cutoffsumhol}
Given a holomorphic regularization ${\mathcal R}: \sigma\mapsto \
\sigma(z)$ on
${\mathcal P}^{*,k}$,  for any $\sigma\in {\mathcal P}^{*,k}$,  the map
 $z\mapsto\cutoffsum_{0}^\infty \sigma(z)$
  is meromorphic with poles of order
 $\leq k+1$  in the discrete
set $\alpha^{-1}\left(\{ -1,0, 1, 2,\cdots \}\right)$  whenever $\sigma(z)$ is a
  holomorphic family with  order $\alpha(z)$ such that ${\rm
      Re}(\alpha^\prime(z))\neq 0$ for any $z$ in $\Omega$.
\end{prop}
 Let $\Omega\subset \C$ be an open
  neighborhood of $0$. Given symbols  $\sigma_1, \cdots ,\sigma_k \in {\mathcal
      P}^{*,0}$, and a holomorphic regularization ${\mathcal R}$ which sends $\sigma_i$
    to $\sigma_i(z)$ with order $\alpha_i(z)$,
    $z\in \Omega$, we build holomorphic perturbations  in the complex multivariable
    $\underline z:=(z_1, \cdots, z_k)\in \Omega^k$ of the symbols $\widetilde \sigma$ introduced in (\mref{eq:sigmatilde}):
\begin{eqnarray*}
 \quad {\widetilde \sigma}(\underline z)&:=&\sigma_1(z_1)\, \overline
 Q\Big(\cdots \sigma_{k-2}(z_{k-2})\, \overline Q\big(\sigma_{k-1}(z_{k-1})\,
\overline Q(\sigma_k(z_k))\big)...\Big).
\end{eqnarray*}
By Theorem \mref{thm:sigmatilde}, these are linear combinations of
log-polyhomogeneous symbols of log-type $j-1$  and
order   $\alpha_1(z_1)+\cdots +\alpha_j(z_j)+j-1$,
$j\in\{1,\ldots ,k\}$.
Applying Proposition \mref{prop:cutoffsumhol} to  each of these symbols provides
information on the pole structure of nested sums of (positively supported)
classical symbols reminiscent of the pole structure of multiple zeta functions \mcite{AET,Go,Zh2}.
\begin{thm}\mlabel{thm:meroChensums}
   Fix symbols  $\sigma_1, \cdots ,\sigma_k \in {\mathcal
      P}^{*,0}$ and  a holomorphic regularization ${\mathcal R}$ which sends $\sigma_i$
    to $\sigma_i(z)$ with order $\alpha_i(z)$.
\begin{enumerate}
\item The map
\begin{eqnarray*}
(z_1, \cdots, z_k)&\mapsto& \altcutoffsum^{\smop{Chen}}_{<}\sigma_1(z_1)\otimes \cdots \otimes
 \sigma_k(z_k)
\end{eqnarray*}
is meromorphic with poles on a countable number of hypersurfaces
$$\sum_{i=1}^j \alpha_i(z_i)\in
-j +\N_0,$$
of multiplicity $j$ varying in $\{1, \cdots, k\}$. Here $\N_0$ stands for the
set of non-negative integers.
\item Let
    $\sigma(z):=\sigma_1(z)\otimes\cdots\otimes \sigma_k(z)$ with $z\in \Omega$. Assume that the orders $\alpha_i(z)$ of the $\sigma_i$'s are nonconstant
    affine with  $\alpha_j^\prime(0)=-q$ for any $j$ in $\{1, \cdots, k\}$
    and some positive real number $q$.
The map $z\mapsto
\altcutoffsum^{\smop{Chen}}_<\sigma(z)$
is meromorphic on $\Omega$ with poles $ z \in (\sum_{i=1}^j \alpha_i(0)+j-\N_0)/(q\, j)$
 of order $\leq j$.
\item If $\mop{Re }(\alpha_1(z_1)+\cdots +\alpha_j(z_j))<-j$ for any $j\in\{1,\ldots
  ,k\}$, the nested sums converge and boil down to ordinary nested  sums (independently
  of the perturbation). Setting $\sigma=\sigma_1\otimes \cdots \otimes \sigma_k$ we have:
\begin{eqnarray*}\cutoffsum^{\smop{Chen}, {\mathcal R}}_{<}
 \sigma&=&
 \lim_{z\to 0} \cutoffsum^{\smop{Chen}}_{<}
 \sigma(z)=
 \sum_<^{\smop{Chen}}\sigma.
\end{eqnarray*}
\end{enumerate}
\end{thm}

\subsection{A twisted holomorphic regularization}

We now take $\calap$ to be a subalgebra  of $ {\mathcal P}^{*,0}$ equipped with
the ordinary product on functions.
Any holomorphic regularization
${\mathcal R}$ on   $\calap$  with parameter space $\Omega\subset \C$  induces  one on the tensor algebra  ${T}(\calap)$:
$$\widetilde {\mathcal R}(\sigma_1\otimes \cdots \otimes \sigma_k)(z_1, \cdots,
z_k):={\mathcal R}(\sigma_1)(z_1)\otimes
\cdots\otimes {\mathcal R}(\sigma_k)(z_k).$$
 It is compatible with the shuffle product
$$\widetilde {\mathcal R}\left((\sigma_1\otimes \cdots \otimes \sigma_k)\, \shu \, (\sigma_{k+1}\otimes
\cdots \otimes \sigma_{k+l})\right)= \widetilde {\mathcal R}\left(\sigma_1\otimes \cdots
\otimes \sigma_k\right)\,\shu  \widetilde {\mathcal R}\left(\sigma_{k+1}\otimes \cdots
\otimes \sigma_{k+l}\right)
$$
for any  $ \sigma_i\in \calap, \, i\in \{1, \cdots, k
+l\}$.
\begin{rk}
Note that $\widetilde {\mathcal R}(\sigma_1\shu
   \sigma_2)(z_1,z_2)
\neq {\mathcal R}(\sigma_1)(z_1)\shu {\mathcal
     R}(\sigma_2)(z_2)$ whereas  $\widetilde {\mathcal R}(\sigma_1\shu
   \sigma_2)(z_1,z_2)
= ({\mathcal R}(\sigma_1)\shu {\mathcal
     R}(\sigma_2))(z_1, z_2)$.
\end{rk}
Let
$$
\delta_k:\C\to \C^{\otimes k}, \quad z\mapsto z\cdot 1^{\otimes k},
$$
be the diagonal map  $\delta: \C\mapsto {T}(\C)$
 and $\delta^*$ the induced map on tensor products of holomorphic symbols
$$\delta_k^\star: {T} \left({\rm Hol}_\Omega\left(\calap\right)\right)\to {\rm Hol}_\Omega\left({T}(\calap)\right),\quad
\sigma \mapsto \sigma\circ \delta_k.
$$
The regularization $\widetilde{\mathcal R}$ induces a one parameter holomorphic regularization:
$$\left(\delta^*\circ \widetilde {\mathcal R}\right)(\sigma_1\otimes \cdots
\otimes \sigma_k)(z)= {\mathcal R}(\sigma_1)(z)\otimes \cdots \otimes {\mathcal R}(\sigma_k)(z)$$
compatible with the shuffle product:
\begin{eqnarray}\mlabel{eq:deltaregshu}&&\left(\delta^*\circ \widetilde {\mathcal R}\right)\left((\sigma_1\otimes \cdots \otimes \sigma_k)\, \shu \, (\sigma_{k+1}\otimes
\cdots \otimes \sigma_{k+l})\right)\nonumber\\
&=&\left(\delta^*\circ \widetilde {\mathcal R}\right)\left(\sigma_1\otimes \cdots
\otimes \sigma_k\right)\,\shu \, \left(\delta^*\circ \widetilde {\mathcal R}\right)\left(\sigma_{k+1}\otimes \cdots
\otimes \sigma_{k+l}\right),
\end{eqnarray}
for any  $ \sigma_i\in \calap, \, i\in \{1, \cdots, k
+l\}$.

Twisting it by
Hoffman's isomorphism in Theorem~\mref{th:hoffman} yields a holomorphic  regularization $\left(\delta^*\circ \tilde {\mathcal R}\right)\ast$ (denoted by ${\mathcal R}^*$ in~\mcite{MP2}) on  ${T}(\calap)$:
$$\left(\delta^*\circ \tilde {\mathcal R}\right)\ast:=\exp\circ \left(\delta^\star
\circ \tilde {\mathcal R}\right)\circ\log,$$
which is compatible with the stuffle product:
\begin{equation}\mlabel{eq:Rstuz}\left(\delta^*\circ \tilde {\mathcal R}\right)^\ast(\sigma\qspr \tau)= \left(\delta^*\circ \tilde {\mathcal R}\right)^\ast(\sigma)\qspr
\left(\delta^*\circ \tilde {\mathcal R}\right)^\ast(\tau),\quad \forall \sigma, \tau \in  {T}\left(\calap\right).
\end{equation}
Consequently, the following regularization
\begin{equation}\mlabel{eq:starR}\widetilde{\mathcal R}^\ast(\sigma_1\otimes \cdots \otimes\sigma_k)(z_1,
\cdots, z_k)=\exp\circ
\tilde {\mathcal R}\circ\log(\sigma_1\otimes \cdots
\otimes\sigma_k)\left(z_{1}, \cdots, z_{k}\right)\end{equation}
is compatible  with stuffle relations after symmetrization in the
complex variables $z_i$
\begin{equation}\mlabel{eq:Rstuzi}\left(\widetilde {\mathcal R}^\ast(\sigma\qspr \tau)\right)_{\rm sym}= \left(\tilde {\mathcal R}^\ast(\sigma)\qspr
\tilde {\mathcal R}^\ast(\tau)\right)_{\rm sym},\quad \forall \sigma, \tau \in  {\mathcal
  T}\left(\calap\right),
\end{equation}
where the subscript {\rm  sym } stands for symmetrization
$$f_{\rm sym}(z_1, \cdots, z_k):= \frac{1}{k!}\sum_{\tau \in
    \Sigma_k}  f(z_{\tau (1)}, \cdots, z_{\tau(k)}),$$  over all the complex
variables $z_1, \cdots, z_{k+l}$
 if $\sigma$ is a tensor of degree $k$ and
$\tau$ a tensor of degree $l$.
Setting $z_1=\cdots=z_{k+l}=z$  in (\mref{eq:Rstuzi}) yields back (\mref{eq:Rstuz})
so that  (\mref{eq:Rstuzi})  can be seen as a polarization of
(\mref{eq:Rstuz}).
Given symbols  $\sigma_1, \cdots ,\sigma_k $ in ${\mathcal
      P}^{*,0}$, and  a holomorphic regularization  ${\mathcal R}: \sigma\mapsto \sigma(z)$,
    sending $\sigma_i$ to $\sigma_i(z)$ with order $\alpha_i(z)$, we are now ready to build a
  map
$$
(z_1, \cdots, z_k) \mapsto \altcutoffsum^{\smop{Chen}}_{<}\widetilde{\mathcal R}^\ast(\sigma_1 \otimes
\cdots \otimes \sigma_k )(z_1, \cdots, z_k),
$$
which, by Theorem \mref{thm:meroChensums},
is meromorphic with poles on a countable number of hypersurfaces
$$\sum_{i=1}^j \alpha_i(z_i)\in
-j +\N_0,$$
with multiplicity $j$ varying in $\{1, \cdots, k\}$.
In particular, if the holomorphic
regularization  ${\mathcal R}$ sends a symbol $\sigma$ to a symbol $\sigma(z)$ with order
$\alpha(z)=\alpha(0)-q\, z$ for some positive real number $q$, the
hypersurfaces of poles are given by
$$\sum_{i=1}^j z_i\in \frac{\sum_{i=1}^j \alpha_i(0)+j-n}{q}, \quad n\in
\N_0, $$
so that hyperplanes of poles containing the origin correspond to $\sum_{i=1}^j
z_i=0$
each of which  with  multiplicity $j$ varying in $\{1, \cdots, k\}$.
\subsection{Meromorphic nested sums of symbols}
Let  ${\rm Mer}_0(\C)$ denote the germ of meromorphic functions in a neighborhood
of zero  in the complex plane  and let ${\rm Hol}_0(\C)$ be the germ of holomorphic
functions at zero. We consider the (Grothendieck closure of the)  tensor algebra
$${T}\left({\rm
    Mer}_0(\C)\right)=\oplus_{k=0}^\infty{T}^k({\rm Mer}_0(\C))$$ over
${\rm Mer}_0(\C)$  and its subalgebra ${T}({\rm Hol}_0(\C)):=
\oplus_{k=0}^\infty{T}^k({\rm Hol}_0(\C))$ where we have set  ${T}^k({\rm Mer}_0(\C)):= \hat \otimes^k {\rm Mer}_0(\C)$, ${T}^k({\rm Hol}_0(\C)):= \hat \otimes^k {\rm Hol}_0(\C)$, and where $\hat
\otimes$ stands for the Grothendieck closure. They come equipped with
the product:
$$(f_1\otimes \cdots \otimes f_k)\bigotimes (f_{k+1}\otimes \cdots \otimes f_{k+l})= f_1\otimes \cdots \otimes f_k\otimes f_{k+1}\otimes \cdots \otimes f_{k+l}.$$
We consider the following  linear extension of ${T}^k({\rm Mer}_0(\C))$ which corresponds to germs at zero of meromorphic maps in
severable variables with linear poles. Let  ${\mathcal L}{\rm Mer}_0(\C^\infty):= \oplus_{k=1}^\infty {\mathcal L}{\rm Mer}_0(\C^k)$ where
$$
{\mathcal L}{\rm Mer}_0(\C^k):= \left\{\prod_{i=1}^m f_i\circ L_i\ \Big | \quad
 f_i\in {\rm Mer}_0(\C), \quad L_i\in \left(\C^k\right)^*\right\}
$$  or equivalently,
$${\mathcal L}{\rm Mer}_0(\C^k):= \left\{(z_1, \cdots, z_k)\mapsto \frac{h(z_1,\cdots, z_k)}{\prod_{L\in \left(\C^k\right)^*}
 \left( L(z_1,\cdots, z_k)\right)^{m_L}}\quad \Big | \quad h \in {\rm Hol
}_0\left(\C^k\right), \quad m_L\in \N \right\}.$$ Setting $m=k$ and $ L_i(z_1, \cdots,
z_k)=z_i $ yields a canonical
injection \begin{eqnarray*}
i:{T}^k({\rm Mer}_0(\C))&\to &{\mathcal L}{\rm Mer}_0(\C^k)\\
f_1\otimes \cdots \otimes f_k&\mapsto&\left( (z_1, \cdots, z_k)\mapsto \prod_{i=1}^k f_i\circ L_i(z_1, \cdots, z_k)\right),
\end{eqnarray*}   and the tensor product on ${T}({\rm Mer}_0(\C))$ extends to
${\mathcal L}{\rm Mer}_0(\C^\infty)$, by
{\small
\begin{eqnarray}\mlabel{eq:bullet}
 && \left((z_1,\cdots, z_k)\mapsto \prod_{i=1}^m f_i\circ L_i (z_1, \cdots,
   z_k)
\right) \bullet \left((z_{1}, \cdots, z_{l}) \mapsto \prod_{j=1}^n f_{m+j}\circ L_{i+j} (z_{1}, \cdots, z_{l}) \right)\\
&=& \left((z_1,\cdots, z_k, \cdots, z_{k+l})\mapsto \prod_{i=1}^m f_i\circ L_i (z_1, \cdots, z_k)\, \prod_{j=1}^n f_{m+j}\circ L_{m+j} (z_{k+1}, \cdots, z_{k+l})\right)\nonumber
\end{eqnarray}
}
which makes it  a graded
algebra.

Specializing to linear forms ${\mathcal L}_k:=\{L\in \left(\C^k\right)^*\ |\  \exists\, J\subset \{1,
\cdots, k\}, \quad  L(z_1,
\cdots, z_k)=\sum_{j\in J}z_j\}$ gives rise to a subalgebra ${\mathcal L} {\mathcal
  M}_0(\C^\infty):= \oplus_{k=1}^\infty {\mathcal L}{\mathcal M}_0(\C^k)\subset {\mathcal
  L} {\rm Mer}_0(\C^\infty)$ defined by:
 $${\mathcal L}{\mathcal M}_0(\C^k):= \left\{(z_1, \cdots, z_k)\mapsto
   \frac{h(z_1,\cdots, z_k)}{\prod_{L\in {\mathcal L}_k}
 \left( L(z_1,\cdots, z_k)\right)^{m_L}}\quad \Big |  \quad h \in {\rm Hol
}_0\left(\C^k\right), \quad m_L\in \N \right\}.$$
For future use, we  consider the   map $\delta^*: {\mathcal L} {\mathcal M}_0(\C^k)\to
{\rm Mer}_0(\C)$ defined by
$$\delta_k^\star: {\mathcal L}{\mathcal M}_0(\C^k)\to {\rm Mer}_0(\C), \quad f\mapsto f\circ \delta_k,
$$
induced by  the diagonal map $\delta: \C\mapsto {T}(\C)$
 previously
defined.

By definition of the twisted regularization  $\tilde{\mathcal
    R}^\ast$, the expressions
 $\cutoffsum^{\smop{Chen}}_<\tilde{\mathcal R}^\ast(\sigma_1\otimes \cdots\otimes
\sigma_k)(z_1, \cdots, z_k)$ are linear combinations of expressions of the
type $\altcutoffsum^{\smop{Chen}}_{<}\tau_1(u_1)\otimes \cdots \otimes
 \tau_l(u_k)$ with  symbols $\tau_j(u_j)$  built from  products of the
 $\sigma_i(z_i)$'s.  It therefore follows from   Theorem
 \ref{thm:meroChensums}, that the functions
$(z_1, \cdots, z_k)\mapsto \cutoffsum^{\smop{Chen}}_<\tilde{\mathcal R}^\ast(\sigma_1\otimes \cdots\otimes
\sigma_k)(z_1, \cdots, z_k)$ lie  in ${\mathcal L}{\mathcal M}_0(\C^k)$.
Since the stuffle relations are satisfied for convergent nested sums,
given  two tensor products $\sigma=
\sigma_1\otimes \cdots\otimes \sigma_k$ and $\tau=\tau_1\otimes \cdots\otimes
\tau_l$ of symbols in ${\mathcal P}$, setting $\sigma_i(z_i):= {\mathcal R}(\sigma_i)(z_i)$,
for ${\rm Re} (z_i)$ sufficiently large we have:
\begin{eqnarray*}
&& \cutoffsum^{\smop{Chen}}_<\left(\tilde{\mathcal R}^\ast\left(\sigma_1\otimes \cdots\otimes
\sigma_k\right)(z_1, \cdots, z_k)\right)\ast \left(\tilde{\mathcal R}^\ast\left(\tau_1\otimes \cdots\otimes
\tau_l\right)(z_{k+1},  \cdots, z_{k+l})\right)\\
&=&\left(\cutoffsum^{\smop{Chen}}_<\tilde{\mathcal R}^\ast(\sigma_1\otimes \cdots\otimes
\sigma_k)(z_1, \cdots, z_k)\right)\, \left( \cutoffsum^{\smop{Chen}}_<\tilde{\mathcal R}^\ast(\tau_1\otimes \cdots\otimes
\tau_l)(z_{k+1}, \cdots,z_{k+l})\right).
\end{eqnarray*}
By analytic continuation (see for example \mcite{GR}, in particular the
Identity Theorem in Chapter 1, Section A, or \mcite{Ho}), this holds as an
identity of meromorphic functions. Since
\begin{eqnarray*}
&&\left(\tilde {\mathcal R}^\ast\left((\sigma_1\otimes \cdots\otimes
\sigma_k)\ast (\tau_1\otimes \cdots\otimes
\tau_l)\right)(z_1, \cdots, z_{k+l})\right)_{\rm sym}\\
&=& \left(\left(\tilde
 {\mathcal R}^\ast(\sigma_1\otimes \cdots\otimes
\sigma_k)\ast \tilde {\mathcal R}^\ast(\tau_1\otimes \cdots\otimes
\tau_)\right) (z_{k+1}, \cdots, z_{k+l})\right)_{\rm sym},
\end{eqnarray*}
 symmetrization in the variables
$z_i$  yields
\begin{eqnarray*}
&& \left(\cutoffsum^{\smop{Chen}}_<\left(\tilde {\mathcal R}^\ast(\sigma_1\otimes \cdots\otimes
\sigma_k)\ast \tilde{\mathcal R}^\ast(\tau_1\otimes \cdots\otimes
\tau_l)\right)(z_{1},  \cdots, z_{k+l})\right)_{\rm sym}\\
&=&\left(\left(\cutoffsum^{\smop{Chen}}_<\tilde{\mathcal R}^\ast(\sigma_1\otimes \cdots\otimes
\sigma_k)(z_1, \cdots, z_k)\right)\, \left( \cutoffsum^{\smop{Chen}}_<\tilde {\mathcal R}^\ast(\tau_1\otimes \cdots\otimes
\tau_l)(z_{k+1}, \cdots,z_{k+l})\right)\right)_{\rm sym}.
\end{eqnarray*}
This can be reformulated as follows.
\begin{thm}\mlabel{thm:meroChensums2} \mcite{MP2} Let $\calap$ be a subalgebra
    of ${\mathcal P}^{*,0}$ and let  ${\mathcal R}$ be a holomorphic regularization which  sends a symbol $\sigma$ to a symbol $\sigma(z)$ with order
$\alpha(z)=\alpha(0)-q\, z$ for some positive real number $q$.
\begin{enumerate}
\item The map
\begin{eqnarray*}
\Psi^{{\mathcal R}}:\left( {T}\left(\calap\right), \qspr      \right)&\to&
\left({\mathcal L}{\mathcal M}_0(\C^\infty), \bullet\right)\\
 \sigma_1\otimes \cdots\otimes \sigma_k &\mapsto & \left((z_1, \cdots, z_k)\mapsto
\cutoffsum^{\smop{Chen}}_<\tilde{\mathcal R}^\ast(\sigma_1\otimes \cdots\otimes
\sigma_k)(z_1, \cdots, z_k)\right),
\end{eqnarray*}
satisfies the following  relation:
\begin{eqnarray*}
\left(\Psi^{\mathcal R}(\sigma\qspr \tau)\right)_{\rm sym}&=&\left(\Psi^{\mathcal R} (\sigma)\,
 \bullet \Psi^{\mathcal R} (\tau)\right)_{\rm sym},
\end{eqnarray*}
which holds as an equality of meromorphic functions in several variables. Here,
as before the subscript sym stands for the symmetrization in the complex
variables $z_i$.
\item After composition with $\delta^*$ this in turn gives rise to a map

\begin{eqnarray}\label{eq:psiR}
\psi^{{\mathcal R}}:\left( {T}\left(\calap\right), \qspr      \right)&\to&
{\rm Mer}_0(\C)\nonumber\\
 \sigma_1\otimes \cdots\otimes \sigma_k &\mapsto & \left(z\mapsto
\delta^*\circ\cutoffsum^{\smop{Chen}}_<\tilde{\mathcal R}^\star(\sigma_1\otimes \cdots\otimes
\sigma_k)(z)\right),
\end{eqnarray}
which is an algebra morphism. In other
words, $\psi^{\mathcal R}$ satisfies the relation:
\begin{eqnarray*}
\psi^{\mathcal R}(\sigma\qspr \tau)&=&\psi^{\mathcal R} (\sigma)\,
 \cdot \psi^{\mathcal R} (\tau),
\end{eqnarray*}
which holds as an equality of meromorphic functions in one variable.
\end{enumerate}
\end{thm}
\subsection{Renormalized nested sums of symbols}
We want to extract finite parts from the meromorphic functions in Theorem \mref{thm:meroChensums2} while
preserving the stuffle relations using a renormalization
procedure. Renormalized evaluators inspired from generalized evaluators used in
  physics  provide a first renormalization procedure.
\subsubsection{Renormalized nested sums via renormalized evaluators}
We call  regularized evaluator at zero   on the germ ${\rm Mer}_0(\C)$ of
meromorphic functions around zero,  any
  linear form  on  ${\rm Mer}_0(\C)$ which extends the evaluation at zero
  ${\rm ev}_0:h\mapsto h(0)$  on
  holomorphic germs at zero. The map ${\rm ev}_0^{\rm reg}$ defined by
$$  {\rm ev}_0^{\rm reg}:= {\rm ev}_0\circ (I-\Pi),$$
where $\Pi:{\rm Mer}_0(\C)\to {\rm Mer}_0(\C)$ as defined in Example
\ref{ex:Pi} corresponds to  the projection onto the pole
part of the Laurent expansion at zero, is such a regularized evaluator at
zero. When we need to specify the complex variable $z$ we also write ${\rm
  ev}_{z=0}^{\rm reg}$.   Following Speer \mcite{S} we introduce renormalized
evaluators which correspond to his generalized evaluators.
\begin{defn} A renormalized evaluator  $\Lambda$ on  a graded subalgebra ${\mathcal
    B}=\oplus_{k=0}^\infty{\mathcal B}_k$ of ${\mathcal L}{\rm
    Mer}_0(\C^\infty)=\oplus_{k=0}^\infty{\mathcal L}{\rm
    Mer}_0(\C^k) $ equipped with the product $\bullet$ introduced in (\mref{eq:bullet}), is
a character  on ${\mathcal B}$ which is compatible with the filtration induced by
the grading and  extends the ordinary evaluation at zero on
holomorphic maps. Equivalently,
\begin{enumerate}
\item Compatibility with the filtration: Let ${\mathcal B}^K:= \oplus_{k=0}^K{\mathcal B}_k$ and   $\Lambda_K:=
  \Lambda_{\vert_{{\mathcal B}^K}}$. Then $\Lambda_{K+1}{}_{{\vert_{{\mathcal
        B}^K}}}=\Lambda_K$.
\item It coincides with the  evaluation map  at zero   on holomorphic maps:
$$\Lambda_{\vert_{{T}\left({\rm Hol}_0(\C)\right)}}={\rm ev}_0.$$
\item It fulfills a multiplicativity property:  $$\Lambda (f\bullet g)= \Lambda(f)\, \Lambda(g),\quad \forall f, g\in {\mathcal
B}.$$
\end{enumerate}
We call the evaluator symmetric if moreover for any $f$ in ${\mathcal B}_k$ and
$\tau$ in $\Sigma_k$,  we have
$$\Lambda(f_\tau)= \Lambda(f),\quad \forall \tau\in \Sigma_k,$$
where we have set
$f_\tau(z_1, \cdots, z_k):=f(z_{\tau(1)}, \cdots, z_{\tau(k)})$.
\end{defn}
\begin{ex} Any regularized evaluator at zero $\lambda$ on ${\rm Mer}_0(\C)$ uniquely extends to a renormalized evaluator $\tilde \lambda$ on the tensor algebra $\left({T}\left({\rm Mer}_0(\C)\right),\otimes\right)$ defined by
$$\tilde \lambda(f_1\otimes \cdots \otimes f_k)=\prod_{i=1}^k \lambda(f_i).$$
\end{ex}
\begin{ex} Any regularized evaluator $\lambda$ on ${\rm Mer}_0(\C)$ extends to
   renormalized evaluators $\Lambda$ and $\Lambda^\prime$  on ${\mathcal L}{\rm Mer}_0(\C^\infty)$
  defined on ${\mathcal L}{\rm Mer}_0(\C^k)$  by
$$\Lambda :=\lambda_{z_{1}}\circ\cdots\circ \lambda_{z_{k}},\quad \Lambda^\prime :=\lambda_{z_{k}}\circ\cdots\circ \lambda_{z_{1}}$$
and to a symmetrized evaluator defined on ${\mathcal L}{\rm Mer}_0(\C^k)$  by
$$\Lambda^{\rm sym} :=\frac{1}{k!}\sum_{\tau\in \Sigma_k}\lambda_{z_{\tau(1)}}\circ\cdots\circ \lambda_{z_{\tau(k)}},$$  where
$\lambda_{z_i}$ stands for the evaluator $\lambda$ implemented in the sole variable
$z_i$, the others being kept fixed.
Their restrictions to $ {T}\left({\rm Mer}_0(\C)\right)$ all  coincide with $\tilde \lambda$.
\end{ex}

\begin{ex}\mlabel{ex:renev} Take $\lambda:= {\rm ev}^{\rm reg}_{0}$, and set with the above notations
$${\rm ev}_0^{\rm ren}:= \Lambda; \quad {\rm ev}_0^{{\rm ren}^\prime}:=
\Lambda^\prime,\quad {\rm ev}_0^{\rm ren,sym}:= \Lambda^{\rm sym},$$ then given a
holomorphic function $h(z_1,z_2)$ in a neighborhood of $0$ and setting
$f(z_1,z_2):= \frac{h(z_1,z_2)}{z_1+z_2}$, we have
$${\rm ev}_0^{\rm ren}\left(f\right)= \partial_1h(0,0); \ {\rm ev}_0^{{\rm ren}^\prime}(f)=
\partial_2h(0,0); \  {\rm ev}_0^{\rm ren,sym}\left(f\right)=\frac{\partial_1h(0,0)+
  \partial_2h(0,0)}{2}=
{\rm ev}_0^{\rm reg}\circ \delta^*\left( f\right),$$
though in general, $${\rm ev}_0^{\rm ren,sym}\neq {\rm ev}_0^{\rm
  reg}\circ \delta^*.$$
\end{ex}

\begin{prop}\mlabel{prop:stufflePsi}
Let $\calap$ be a subalgebra
    of ${\mathcal P}^{*,0}$ and let  ${\mathcal R}$ be a holomorphic regularization which  sends a symbol $f$ to a symbol $\sigma(z)$ with order
$\alpha(z)=\alpha(0)-q\, z$ for some positive real number $q$.
Let ${\mathcal E}$ be a {\rm symmetrized} renormalized evaluator on ${\mathcal L}{\mathcal M}_0$.  The map
\begin{eqnarray*}
\Psi^{{\mathcal R},{\mathcal E}}:\left( {T}\left(\calap\right), \qspr\right)&\to& \C\\
 \sigma_1\otimes \cdots\otimes \sigma_k &\mapsto &
{\mathcal E}\circ\Psi^{{\mathcal R}} (\sigma_1\otimes \cdots \otimes \sigma_k)
\end{eqnarray*}
defines a character. In other
words, $\Psi^{{\mathcal R}, {\mathcal E}}$ satisfies the stuffle relation:
\begin{eqnarray*}
\Psi^{{\mathcal R}, {\mathcal E}}(\sigma\qspr \tau)&=&\Psi^{{\mathcal R}, {\mathcal E}} (\sigma)\,
 \cdot\, \Psi^{{\mathcal R}, {\mathcal E}} (\tau).
\end{eqnarray*}
\end{prop}
\begin{rk}
{\rm Here, we use the
fact that for a symmetrized evaluator $\Lambda$ we have $\Lambda(f)= \Lambda(f_{\rm sym})$ where as before the subscript ``sym'' stands for the symmetrization in
the complex variables $z_i$.
}
\end{rk}

This proposition gives rise to renormalized nested sums of symbols
$$\cutoffsum_{<}^{\rm Chen,{\mathcal R}, {\mathcal E}} \sigma_1\otimes \cdots\otimes \sigma_k:= \Psi^{{\mathcal R}, {\mathcal E}}(\sigma_1\otimes \cdots\otimes \sigma_k)$$
which obey stuffle relations:
$$\cutoffsum_{<}^{\rm Chen,{\mathcal R}, {\mathcal E}} (\sigma\qspr \tau)= \left(\cutoffsum_{<}^{\rm Chen,{\mathcal R}, {\mathcal E}} \sigma\right)\, \left(\cutoffsum_{<}^{\rm Chen,{\mathcal R}, {\mathcal E}} \tau\right).$$

\subsubsection{Renormalized nested sums via algebraic Birkhoff decomposition}
On the other hand, the tensor algebra  ${T}(\calap)$ can be equipped with the deconcatenation coproduct:
$$\Delta \left(\sigma_1\otimes \cdots \otimes \sigma_k\right):= \sum_{j=0}^k
\left(\sigma_1\otimes\cdots \otimes \sigma_j\right) \bigotimes\left(
\sigma_{j+1}\otimes \cdots \otimes \sigma_k\right)$$  which then inherits  a structure of connected
graded commutative Hopf algebra \mcite{Ho2}.
Using the
convolution product $\star$ associated with the product and coproduct on
${T}(\calap)$ and since  ${\rm Mer}_0(\C)$ embeds into the Rota-Baxter algebra $\CC[\vep^{-1},\vep]]$ we can  implement an algebraic Birkhoff decomposition as in (\ref{eq:ABF})
to the map $\psi^{{\mathcal R}}$ in Eq.~(\mref{eq:psiR}):
$$\psi^{{\mathcal R}}= \left(\psi^{{\mathcal R}}_-
\right)^{\star (-1)} \star \psi^{{\mathcal R}}_+$$
associated with the minimal
substraction scheme  to build
characters $$ \psi^{{\mathcal R}}_+(0):\left( {T}\left(\calap\right),
  \qspr\right)\to \C.$$

\begin{prop}\mlabel{prop:stufflepsi}\mcite{MP2}
Let $\calap$ be a subalgebra
    of ${\mathcal P}^{*,0}$ and let  ${\mathcal R}$ be a holomorphic regularization which  sends a symbol $\sigma$ to a symbol $\sigma(z)$ with order
$\alpha(z)=\alpha(0)-q\, z$ for some positive real number $q$.  The map
\begin{eqnarray*}
\psi^{{\mathcal R},{\rm Birk}}:\left( {T}\left(\calap\right), \qspr\right)&\to& \C\\
 \sigma_1\otimes \cdots\otimes \sigma_k &\mapsto &
\psi^{{\mathcal R}}_+(0) (\sigma_1\otimes \cdots \otimes \sigma_k)
\end{eqnarray*}
defines a character
\begin{eqnarray*}
\psi^{{\mathcal R}, {\rm Birk}}(\sigma \qspr
\tau)&=&\psi^{{\mathcal R},
 {\rm Birk}} (\sigma)\,
 \cdot\, \psi^{{\mathcal R}, {\rm Birk}} (\tau).
\end{eqnarray*}
\end{prop}
The map yields an alternative set of  renormalized nested sums of symbols
$$\cutoffsum_{<}^{{\rm Chen},{\mathcal R}, {\rm Birk}} \sigma_1\otimes \cdots\otimes \sigma_k:= \psi^{{\mathcal R}, {\rm Birk}}(\sigma_1\otimes \cdots\otimes \sigma_k)$$
which obey stuffle relations:
$$\cutoffsum_{<}^{{\rm Chen},{\mathcal R}, {\rm Birk}} (\sigma \qspr \tau)= \left(\cutoffsum_{<}^{{\rm Chen},{\mathcal R}, {\rm Birk}} \sigma\right)\, \left(\cutoffsum_{<}^{\rm Chen,{\mathcal R}, {\rm Birk}} \tau\right).$$

\subsection{Renormalized (Hurwitz) multiple zeta values at non-positive integers}
\subsubsection{An algebra of symbols}
Since we consider both zeta and Hurwitz zeta functions, let us first observe
that for any non-negative number $v$ and any $\sigma $ in ${\mathcal P}^{*, k}$,
the map $\xi\mapsto t_v^*\sigma(\xi):= \sigma(\xi+v)$
 defines a symbol in ${\mathcal P}^{*, k}$.

Let  $\widetilde\calap$ be the subalgebra of ${\mathcal P}^{*,0}$ generated by
the continuous functions with support inside the interval $(0,1)$ and the set
$$\{\sigma\in {\mathcal P}^{*,0}\,\big|\, \exists v\in [0, +\infty), \exists s\in\C,\, \sigma(\xi)=
(\xi+v)^{-s}\hbox{ when }\xi\ge 1\}.$$ Consider the ideal ${\mathcal
N}$ of $\widetilde\calap$ of continuous functions with support
inside the interval $(0,1)$. The quotient algebra ${\mathcal
  A}=\widetilde \calap/{\mathcal N}$ is then generated by  elements
$\sigma_{s,v} \in {\mathcal P}^{*,0}$ with  $\sigma_{s,v}(\xi)=(\xi+v)^{-s}$
for $|\xi|\ge 1$. For any $v\in\R_+$ the subspace $\calap_v$ of
${\mathcal A}$ generated by $\{\sigma_{s,v}\ |\,s\in\C\}$ is a
subalgebra of $\calap$. We equip $\calap_v$ with the following holomorphic
regularization on an open neighborhood $\Omega$  of $0$ in $\C$:
\begin{eqnarray*}
{\mathcal R}: \calap_v&\to &{\rm Hol}_\Omega\left({\mathcal A }_v\right)\\
\sigma_{s,v} &\mapsto &\left(z\mapsto  (1-\chi)\, \sigma_{s,v}+ \chi\  \sigma_{s+z, v}\right)
\end{eqnarray*}
where $\chi$ is any smooth cut-off function which is identically one outside
the unit ball and vanishes in a small neighborhood of $0$.

Let ${\mathcal W}$ be the $\C$-vector space freely spanned symbols
indexed by sequences $(u_1,\ldots,u_k)$ of real numbers. In other
words, ${\mathcal W}$ is ${T}(W)$ where $W=\oplus_{u\in \RR} \RR
x_u$ where we identify $x_u$ with  $u$ for simplicity and set
$x_u\cdot x_v = x_{u+v}, u,v\in \RR$. We then define the stuffle
product on ${\mathcal W}$ as usual in Eq.~(\mref{eq:quasi}) or
Eq.~(\mref{eq:mapshr}) with $\lambda=1$. The map
$$ \sigma: {\mathcal W} \to {T}({\mathcal A}_v),\quad u=(u_1,\cdots,u_k) \mapsto\sigma_{u;v}:=\sigma_{(u_1,\ldots ,u_k; \, v)}:=\sigma_{u_1;
  v}\otimes\cdots\otimes\sigma_{u_k; v}$$
induces a stuffle product on ${\mathcal
  T}(\calap_v)$:
$$\sigma_{u; v}\qspr \sigma_{u^\prime; v}:= \sigma_{u\qspr
  u^\prime; v}.$$

As before, we twist the regularization $\widetilde {\mathcal R}$ induced by ${\mathcal
  R}$ on ${T}(\calap_v)$ by the Hoffman isomorphism (\mref{eq:his}) to build a twisted
holomorphic regularization $\widetilde {\mathcal R}^\ast$ in several variables which satisfies
$$ \left(\widetilde {\mathcal R }^\ast(\sigma_{u; v})\qspr \widetilde {\mathcal R
  }^\ast(\sigma_{u^\prime; v})\right)_{\rm sym}=  \left( \widetilde {\mathcal R }^\ast(\sigma_{u\qspr
  u^\prime; v})\right)_{\rm sym}$$
and a twisted
holomorphic regularization $\delta^*\circ \widetilde {\mathcal R}^\ast$ in one variable compatible with the stuffle
product:
$$ \left(\delta^*\circ\widetilde {\mathcal R }^\ast(\sigma_{u; v})\right)\qspr\left( \delta^*\circ \widetilde {\mathcal R }^*(\sigma_{u^\prime; v})\right)=  \delta^*\circ\widetilde {\mathcal R }^\ast(\sigma_{u\qspr
  u^\prime; v}).$$

\subsubsection{Multiple zeta values renormalized via renormalized evaluators}
Let $\Omega $ be an open neighborhood of $0$ in $\C$ and let ${\mathcal R}:\sigma\mapsto
\{\sigma(z)\}_{z\in \Omega}$ be  the   holomorphic regularization procedure on $\widetilde{\mathcal
  A}$  previously introduced. The   multiple Hurwitz zeta functions   defined by:
\begin{eqnarray*}
 \zeta (s_1, \ldots, s_k; \, v_1, \ldots, v_k)&:=&{\Psi}^{{\mathcal R} } (\sigma_{s_1, v_1}\otimes \cdots \otimes\sigma_{s_k, v_k})
\end{eqnarray*}
are
meromorphic in all variables with poles\footnote{When $k=2$ and  $v_1=\cdots=
v_l=v$ a more refined analysis actually shows that  for some any negative real number $v$, poles
actually only arise for  $s_1=-1$ or
  $s_1+s_2\in \{-2,-1, 0, 2,4, 6, \cdots\}$.} on a countable family of hyperplanes $s_1+\cdots +s_j\in
]-\infty, j]\cap \Z$, $j$ varying from $1$ to $k$. When $v_1=\cdots=v_k=v$, we
set $$ \zeta(s_1, \ldots, s_k; \, v):=\zeta (s_1, \ldots, s_k; \, v_1, \ldots,
v_k)$$
in which case they satisfy the following
relations:
$$
\left(\zeta^{{\mathcal E} } (u\qspr u^\prime; v)\right)_{\rm sym} =
 \left( \zeta^{ {\mathcal E } } (u;v) \, \zeta^{ {\mathcal E} }(u^\prime;
   v)\right)_{\rm sym}.
$$

The renormalized multiple Hurwitz zeta values derived from a {\it symmetrized} renormalized
 evaluator ${\mathcal E}$  on ${\mathcal L}{\mathcal M}_0(\C^\infty)$:
\begin{eqnarray*}
 \zeta^{ {\mathcal E}}(s_1, \ldots, s_k; \, v_1, \ldots, v_k)&:=&{\Psi}^{{\mathcal
    R}, {\mathcal E}}(\sigma_{s_1, v_1}\otimes \cdots \otimes\sigma_{s_k, v_k})
\end{eqnarray*}
denoted by
$\zeta^{{\mathcal R}, {\mathcal E} }(s_1, \ldots, s_k; \, v )$
 when $v_1=\cdots=v_k=v$,
satisfy stuffle relations in that case:
$$\zeta^{{\mathcal E} } (u\qspr u^\prime; v) =
  \zeta^{ {\mathcal E } } (u;v) \, \zeta^{ {\mathcal E} }(u^\prime; v).
$$
Let us compute  renormalized values in   the case $k=2$ using  a renormalized
evaluator. For any $a\in\R$ and $m\in\N-\{0\}$ we
introduce  the notation:
$$
[a]_j:=a(a-1)\cdots(a-j+1).
$$
We extend this
 to $j=0$ and $j=-1$ by setting:
%\begin{equation}
$[a]_0:=1,[a]_{-1}:=\frac{1}{a+1}.$
%\end{equation}\\
 Combining Definition (\mref{eq:starR})$$\widetilde {\mathcal R}^* (\sigma_1\otimes \sigma_2)(z_1,z_2)=
 \sigma_1(z_1)\otimes \sigma_2(z_2)-\frac{1}{2}( \sigma_1\bullet
 \sigma_2)(z_1)+ \frac{1}{2} \sigma_1(z_1)\bullet
 \sigma_2(z_2)$$ applied to  the regularization  $${\mathcal R}(\sigma_i)(z)(x)= (x+v)^{-s_i-z}\quad {\rm of}\quad {\rm
order}\quad  \alpha_i(z)= -s_i-z_i,$$    with the Euler-MacLaurin formula (\ref{eq:EML}),
and following  \mcite{MP2} (see   the proof of Theorem 9),  we compute
 \begin{eqnarray*}
 \zeta(s_1,  s_2; v)(z_1,z_2)&=& \Psi^{\mathcal
  R}(\sigma_{s_1,v}\otimes \sigma_{s_2,v})(z_1,z_2)\nonumber\\
&=& \altcutoffsum_{<}^{\smop{Chen}}\sigma_1(z_1)\otimes\sigma_2(z_2)+\frac{1}{2} \sigma_1(z_1)\, \sigma_2(z_2)-\frac{1}{2}  \, (\sigma_1\,
  \sigma_2)(z_1)\nonumber \\
%&=&
%\sum_{j=0}^{2J_2} B_j\,\frac{ [\alpha_2(z_2)]_{j-1}}{j!}\left(
%\zeta(-\alpha_1(z_1)-\alpha_2(z_2)+j-1;\, v)  -\zeta(-\alpha_1(z_1);\,
%v)\right)\nonumber\\
%&+&\frac{1}{2} \zeta(s_1+s_2+z_1+z_2;\, v)-\frac{1}{2}  \, \zeta(s_1+s_2+z_1;\,% v)\\
%&\hbox to 3mm{}&+
%\frac{ [\alpha_2(z_2)]_ {2J_2+1}}{(2J_2+1)!}\cutoffsum_0^\infty\left(
%(n+v)^{\alpha_1(z_1)}\,  \int_1^n \overline{B_{2J_l+1}}
%(y)\,(y+v)^{\alpha_2(z_2)-2J_2-1} \, dy\right)\nonumber\\
&=&
\sum_{j=0}^{2J_2} B_j\,\frac{ [-s_2-z_2]_{j-1}}{j!}\left(
\zeta(s_1+s_2+z_1+z_2+j-1;\, v)  -\zeta(s_1+z_1;\, v)\right)\\
&&+\frac{1}{2} \zeta(s_1+s_2+z_1+z_2;\, v)-\frac{1}{2}  \, \zeta(s_1+s_2+z_1;\, v)\\
&\hbox to 3mm{}&+
\frac{ [-s_2-z_2]_ {2J_2+1}}{(2J_2+1)!}\cutoffsum_0^\infty\left(
(n+v)^{-s_1-z_1}\,  \int_1^n \overline{B_{2J_l+1}}
(y)\,(y+v)^{-s_2-z_2-2J_2-1} \, dy\right).\nonumber\\
\end{eqnarray*}
Hence, for non-positive integers $s_1=-a_1, s_2=-a_2$ and $2J_2= a_1+a_2+2$  we have:
 \begin{eqnarray}\mlabel{eq:useful}
 \zeta(-a_1,  -a_2; v)(z_1,z_2)&=&
\sum_{j=0}^{a_1+a_2+2} B_j\,\frac{ [a_2-z_2]_{j-1}}{j!}\left(
\zeta(-a_1-a_2+z_1+z_2+j-1;\, v)  -\zeta(-a_1+z_1;\, v)\right)\nonumber\\
&&+\frac{1}{2} \zeta(-a_1-a_2+z_1+z_2;\, v)-\frac{1}{2}  \, \zeta(-a_1-a_2+z_1;\, v)\\
&&+
\frac{ [a_2-z_2]_{a_2+2}}{(a_2+2)!}\cutoffsum_0^\infty\left(
(n+v)^{a_1-z_1}\,  \int_1^n \overline{B_{a_1+a_2+3}}
(y)\,(y+v)^{-2} \, dy\right).\nonumber
\end{eqnarray}
The last line on the r.h.s. is a holomorphic expression at zero on which all
renormalized evaluators at zero  vanish. The second line on the r.h.s is a
linear combination of ordinary zeta functions at negative integers which are
holomorphic at zero. Any
evaluator $\Lambda$ at zero vanishes on these terms; indeed, we have $\Lambda\left( \zeta(-a_1-a_2+z_1+z_2;\,
  v)- \, \zeta(-a_1-a_2+z_1;\, v)\right)= \zeta(-a_1-a_2;\,
  v)- \, \zeta(-a_1-a_2;\, v)=0$.  Only when
evaluated on
the expression on the first line of the r.h.s can various evaluators differ.

We want to implement the symmetrized evaluator at zero  $${\rm ev}_0^{\rm ren,sym}:= \frac{1}{2}\left({\rm ev}^{\rm reg}_{z_2=0}\circ  {\rm ev}^{\rm
  reg}_{z_1=0}+{\rm ev}^{\rm reg}_{z_1=0}\circ  {\rm ev}^{\rm
  reg}_{z_2=0} \right)$$
  introduced in Example  \mref{ex:renev}.  We first compute
\begin{eqnarray*}&&{\rm ev}^{\rm reg}_{z_1=0}\left(  {\rm
   ev}^{\rm reg}_{z_2=0}\left( \zeta(-a_1,-a_2; v)(z_1,z_2)\right)\right)\\
 &=&{\rm ev}^{\rm reg}_{z_1=0}\left(  {\rm
   ev}^{\rm reg}_{z_2=0}  \left(
\sum_{j=0}^{a_1+a_2+2} B_j\,\frac{ [a_2-z_2]_{j-1}}{j!}\left(
\zeta(-a_1-a_2+z_1+z_2+j-1;\, v)  -\zeta(-a_1+z_1;\, v\right)\right)\right)\\
&=&  {\rm
   ev}^{\rm reg}_{z_1=0} \left(
\sum_{j=0}^{a_2+1} B_j\, \frac{ [a_2]_{j-1}}{j!}\, \left(
\zeta(-a_1-a_2+z_1+j-1;\, v)  -\zeta(-a_1+z_1;\, v)\right)\right)       \\
&=& \frac{1}{a_2+1}\sum_{j=0}^{a_2+1} B_j\,{a_2+1\choose j}\, \left(
\zeta(-a_1-a_2+j-1;\, v)  -\zeta(-a_1;\, v)\right).\end{eqnarray*}
When $v=0$ this yields:
\begin{eqnarray}\mlabel{eq:ev12}
 &&{\rm ev}^{\rm reg}_{z_1=0}\left(  {\rm
   ev}^{\rm reg}_{z_2=0}\left( \zeta(-a_1,-a_2)(z_1,z_2)\right)\right)
:= {\rm ev}^{\rm reg}_{z_1=0}\left(  {\rm
   ev}^{\rm reg}_{z_2=0}
 \left(\zeta(-a_1,-a_2)(z_1,z_2;0)\right)\right)\nonumber
\\
&=& \frac{1}{a_2+1}\sum_{j=0}^{a_2+1} B_j\,{a_2+1\choose j}\, \left(
  -\frac{B_{a_1+a_2-j+2}}{a_1+a_2-j+2}
  +\frac{B_{a_1+1}}{a_1+1}\right).
\end{eqnarray}
We next compute
{\allowdisplaybreaks
\begin{eqnarray*}
 &&{\rm ev}^{\rm reg}_{z_2=0}\left(  {\rm
   ev}^{\rm reg}_{z_1=0}\left( \zeta(-a_1,-a_2; v)(z_1,z_2)\right)\right)\\
 &=&
{\rm ev}^{\rm reg}_{z_2=0}\left(  {\rm
   ev}^{\rm reg}_{z_1=0} \left(\frac{B_0 }{a_2-z_2+1}\left(
\zeta(-a_1-a_2+z_1+z_2-1;\, v)  -\zeta(-a_1+z_1;\, v)\right)\right)\right)\\
&&+{\rm ev}^{\rm reg}_{z_2=0}\left(  {\rm
   ev}^{\rm reg}_{z_1=0}  \left(
\sum_{j=1}^{a_1+1} B_j\,\frac{ [a_2-z_2]_{j-1}}{j!}\left(
\zeta(-a_1-a_2+z_1+z_2+j-1;\, v)  -\zeta(-a_1+z_1;\, v)\right)\right)\right)\\
&&+{\rm ev}^{\rm reg}_{z_2=0}\left(  {\rm
   ev}^{\rm reg}_{z_1=0}  \left(
\sum_{j=a_1+2}^{a_1+a_2+2} B_j\,\frac{ [a_2-z_2]_{j-1}}{j!}\left(
\zeta(-a_1-a_2+z_1+z_2+j-1;\, v)  -\zeta(-a_1+z_1;\, v)\right)\right)\right)\\
&=&{\rm ev}^{\rm reg}_{z_2=0}   \left(\frac{B_0 }{a_2+1}\left(
\zeta(-a_1-a_2+z_2-1;\, v)  -\zeta(-a_1;\, v)\right)\right)\\
&&+{\rm ev}^{\rm reg}_{z_2=0}\left(
\sum_{j=1}^{a_1+1} B_j\,\frac{ [a_2-z_2]_{j-1}}{j!}\left(
\zeta(-a_1-a_2+z_2+j-1;\, v)  -\zeta(-a_1;\, v)\right)\right)\\
&&+
\sum_{j=1}^{a_2+1} B_{j+a_1+1}\,\partial_{z_2}\left(\frac{
    [a_2-z_2]_{j+a_1}}{(j+a_1+1)!}\right)_{\vert_{z_2=0}}\,{\rm Res}_{z_2=0}\left(
\zeta(-a_2+z_2+j;\, v) \right)\\
&=&
\frac{1}{a_2+1}\sum_{j=0}^{a_2+1} B_j\,{a_2+1\choose j}\, \left(
\zeta(-a_1-a_2+j-1;\, v)  -\zeta(-a_1;\, v)\right)\\
&&+(-1)^{a_1+1} a_1!a_2!\,\frac{ B_{a_1+a_2+2}}{(a_1+a_2+2)!},
\end{eqnarray*}
}
since the only contribution to the residue comes from the term $j=a_1+a_2+2$.
Since $\zeta(-a):= \zeta(-a;0)= -\frac{B_{a+1} }{a+1}$, this combined with
(\mref{eq:sumBi})  applied to $k= a_2+1$ yields
\begin{eqnarray}\mlabel{eq:ev21}
&&{\rm ev}^{\rm reg}_{z_2=0}\left(  {\rm
   ev}^{\rm reg}_{z_1=0}\left( \zeta(-a_1,-a_2)(z_1,z_2)\right)\right)
:= {\rm ev}^{\rm reg}_{z_2=0}\left(  {\rm
   ev}^{\rm reg}_{z_1=0}\left( \zeta(-a_1,-a_2; 0)(z_1,z_2)\right)\right)\nonumber\\
&=&-\frac{1}{a_2+1}\sum_{j=0}^{a_2+1} B_j\,{a_2+1\choose j}\,
  \frac{B_{a_1+a_2-j+2}}{a_1+a_2-j+2}  +\frac{B_{a_1+1}}{a_1+1}\,\frac{B_{a_2+1}}{a_2+1}\\
&&+(-1)^{a_1+1} a_1!a_2!\,\frac{ B_{a_1+a_2+2}}{(a_1+a_2+2)!}.\nonumber
\end{eqnarray}
Combining (\mref{eq:ev12}) and (\mref{eq:ev21}) yields
\begin{eqnarray}\mlabel{eq:zetaev}
 &&\zeta^{{\rm ev}} (-a_1,-a_2) :={\rm ev}_0^{\rm ren,sym} \left( \zeta(-a_1,-a_2)(z_1,z_2)\right) \nonumber\\
&=&-\frac{1}{a_2+1}\sum_{j=0}^{a_2+1} B_j\,{a_2+1\choose j}\,
  \frac{B_{a_1+a_2-j+2}}{a_1+a_2-j+2} \\
&&+\frac{B_{a_1+1}}{a_1+1}\,\frac{B_{a_2+1}}{a_2+1} +(-1)^{a_1+1} a_1!a_2!\,\frac{ B_{a_1+a_2+2}}{2\, (a_1+a_2+2)!}. \nonumber
\end{eqnarray}
Renormalized multiple zeta values of depth $2$  at non-positive arguments obtained this way,
are rational linear combinations  of Bernoulli numbers, and hence rational
numbers.
More generally, an inductive procedure on $k$ carried out in the same spirit
as the proof of  Theorem 10 in \mcite{MP2} shows that the renormalized multiple zeta values
$\zeta^{ {\mathcal E}}(s_1, \cdots, s_k;v)$ are rational values at non-positive
integer arguments $s_1, \cdots, s_k$  whenever $v$ is rational.
\subsubsection{Multiple zeta values renormalized via Birkhoff decomposition}

 The renormalized multiple Hurwitz zeta values derived from a Birkhoff decomposition:
\begin{eqnarray*}
 \zeta^{ {\rm Birk}}(s_1, \ldots, s_k; \, v_1, \ldots, v_k)&:=&{\Psi}^{{\mathcal
    R}, {\rm Birk}}(\sigma_{s_1, v_1}\otimes \cdots \otimes\sigma_{s_k, v_k})
\end{eqnarray*}
denoted by
$\zeta^{ {\rm Birk} }(s_1, \ldots, s_k; \, v )$ when $v_1=\cdots=v_k=v$,
satisfy stuffle relations
$$
\zeta^{{\rm Birk} } (u\qspr u^\prime; v) =
  \zeta^{ {\rm Birk } } (u;v) \, \zeta^{ {\rm Birk} }(u^\prime; v).
$$

A striking holomorphy property arises at non-positive integer arguments
\mcite{MP2} after implementing the diagonal map $\delta$.
\begin{prop} At non-positive integer arguments $s_i$ and for a rational
  parameter $v$, the map
  $$ z\mapsto  {\psi}^{{\mathcal R} }\left( \sigma_{s_1, v}\otimes \cdots \otimes\sigma_{s_k, v}\right)(z) $$
defined in  (\ref{eq:psiR}) is  holomorphic at zero.
\end{prop}
Consequently,
$$\zeta^{ {\rm Birk} }(s_1, \ldots, s_k; \, v )=\lim_{z \to 0} {\psi}^{{\mathcal
    R}}(\sigma_{s_1, v}\otimes \cdots \otimes\sigma_{s_k, v}).$$
Let us compute double zeta values at non-positive integer arguments using
Birkhoff decomposition.
Setting $z_1=z_2=z$ in (\mref{eq:useful}) leads to
\begin{eqnarray*}
 \zeta(-a_1,  -a_2; v)(z)&=&
\sum_{j=0}^{a_2+1} B_j\,\frac{ [a_2-z]_{j-1}}{j!}\left(
\zeta(-a_1-a_2+2z+j-1;\, v)  -\zeta(-a_1+z;\, v)\right)\nonumber\\
&&+
\frac{ [a_2-z]_{a_2+2}}{(a_2+2)!}\cutoffsum_0^\infty\left(
(n+v)^{a_1-z}\,  \int_1^n \overline{B_{a_1+a_2+3}}
(y)\,(y+v)^{-2} \, dy\right).
\end{eqnarray*}
Evaluating this expression at $z=0$ in a similar manner to the previous
computation, yields:
\begin{eqnarray*}
 \zeta^{\rm Birk}(-a_1,  -a_2; v)&=&\lim_{z\to 0}\left(
\sum_{j=0}^{a_2+1} B_j\,\frac{ [a_2-z]_{j-1}}{j!}\left(
\zeta(-a_1-a_2+2z+j-1;\, v)  -\zeta(-a_1+z;\, v)\right)\right)\nonumber\\
&=&
\sum_{j=0}^{a_2+1} B_j\,\frac{ [a_2]_{j-1}}{j!}\left(
\zeta(-a_1-a_2+j-1;\, v)  -\zeta(-a_1;\, v)\right)\\
 && + (-1)^{a_1+1} a_1!a_2!\,\frac{ B_{a_1+a_2+2}}{2\, (a_1+a_2+2)!}.
\end{eqnarray*}
When $v=0$ this yields  \mcite{MP2}:
\begin{eqnarray*}
 \zeta^{{\rm Birk}}(-a_1,  -a_2)
&:=&  \zeta^{{\rm Birk}}(-a_1,  -a_2;0)\\
&=&-\frac{1}{a_2+1}\sum_{j=0}^{a_2+1} B_j\,{a_2+1\choose j}\,
  \frac{B_{a_1+a_2-j+2}}{a_1+a_2-j+2}  \nonumber\\
&+&\frac{B_{a_1+1}}{a_1+1}\,\frac{B_{a_2+1}}{a_2+1}+ (-1)^{a_1+1} a_1!a_2!\,\frac{ B_{a_1+a_2+2}}{2\, (a_1+a_2+2)!}
\end{eqnarray*}
which coincides with (\mref{eq:zetaev}).

Thus, renormalized double zeta values at non-positive integers obtained by two
different methods -- using
the symmetrized renormalized evaluator ${\rm ev}_0^{\rm ren, sym}$ or a Birkhoff decomposition--
coincide.

Formula (\mref{eq:zetaev}) yields the following table of values
$\zeta(-a_1, -a_2)$ for $a_1,a_2\in \{0, \ldots, 6\}$ derived in \mcite{MP2}:

\begin{equation}\hskip -8mm
\begin{array}{c|c|c|c|c|c|c|c|}
\zeta(-a,-b)&a=0 &a=1&a=2&a=3&a=4&a=5&a=6\\ \hline
&&&&&&&\\
b=0 &\frac{3}{8}&\frac{1}{12}&\frac{7}{720}&-\frac{1}{120}
&-\frac{11}{2\, 520}&\frac{1}{252}&\frac{1}{224}\\
&&&&&&&\\\hline
&&&&&&&\\
b=1 &\frac{1}{24}&\frac{1}{288}&-\frac{1}{240}&-\frac{19}{10\, 080}
&\frac{1}{504}&\frac{41}{20\, 160}&-\frac{1}{480}\\
&&&&&&&\\ \hline
&&&&&&&\\
b=2&-\frac{7}{720}&-\frac{1}{240}&0&\frac{1}{504}&\frac{113}{151\, 200}&-\frac{1}{480}&-\frac{307}{166\,320}\\
&&&&&&&\\ \hline
&&&&&&&\\
b=3&-\frac{1}{240}&\frac{1}{840}&\frac{1}{504}&\frac{1}{28\,800}&-\frac{1}{480}&-\frac{281}{332\,640}&\frac{1}{264}\\
&&&&&&&\\ \hline
&&&&&&&\\
b=4&\frac{11}{2\,520}&\frac{1}{504}&-\frac{113}{151\,200}&-\frac{1}{480}&0&\frac{1}{264}&\frac{117\,977}{75\,675\,600}\\
&&&&&&&\\ \hline
&&&&&&&\\
b=5&\frac{1}{504}&-\frac{103}{60\,480}&-\frac{1}{480}&\frac{1}{1232}&\frac{1}{264}&\frac{1}{127\,008}&-\frac{691}{65\,520}\\
&&&&&&&\\ \hline
&&&&&&&\\
b=6&-\frac{1}{224}&-\frac{1}{480}&\frac{307}{166\,320}&\frac{1}{264}&-\frac{117\,977}{75\,675\,600}&-\frac{691}{65\,520}&0\\
&&&&&&&\\ \hline
\end{array}
\end{equation}
This table of values  differs from the one derived  in \mcite{GZ} (see Table~(\mref{eq:gzt})) with which it
however matches for arguments $(a,b)$
 with  $a+b$
odd and $b\not =0$ and for diagonal arguments  $(-a,-a)$.

%=======================================================================================
%========================================================================================
%========================================================================================
%\addcontentsline{toc}{section}{\numberline {}References}

\end{document}